\documentstyle[11pt,leqno,amscd,amssymb,pstricks,latexsym,amsbsy,xypic,mathrsfs,verbatim]{amsart}

\input xy
\xyoption{all}

\newcommand{\scr}[1]{\mathscr #1}
\newtheorem{Thm}{Theorem}[section]
\newtheorem{Lem}[Thm]{Lemma}
\newtheorem{Cor}[Thm]{Corollary}
\newtheorem{Prop}[Thm]{Proposition}

 \theoremstyle{definition}

\newtheorem{Defn}[Thm]{Definition}

\newtheorem{Rem}[Thm]{Remark}
\newtheorem{Rems}[Thm]{Remarks}

\numberwithin{equation}{section}

\def\nh{{\rm nh}}
\def\Hom{{\rm Hom}}
\def\End{{\rm End}}
\def\Ext{{\rm Ext}}
\def\Aut{{\rm Aut}}

\newcommand{\ul}{\underline}
\newcommand{\co}{{\mathcal O}}
\newcommand{\vx}{\vec{x}}
\newcommand{\vu}{\vec{u}}
\newcommand{\vv}{\vec{v}}
\newcommand{\vc}{\vec{c}}
\newcommand{\vw}{\vec{\omega}}

\def\coh{\mbox{\rm coh-}}
\def\cohz{{\rm coh}_0\mbox{-}}
\def\vect{\mbox{\rm vect-}}
\def\tor{{\rm coh}_0\mbox{-}}
\def\rk{{\rm rk}}

\def\az{\alpha}
\def\bz{\beta}
\def\gz{\gamma}
\def\thz{\theta}
\def\dz{\delta}

\def\lz{\lambda}
\def\sz{\sigma}

\def\dz{\delta}
\def\vphi{\varphi}

\def\Thz{\Theta}
\def\Dz{\Delta}

\def\bfH{{\bf H}}

\def\bfaz{\boldsymbol{\alpha}}  \def\bfbz{\boldsymbol{\beta}}
\def\bfgz{\boldsymbol{\gamma}}  \def\bfxi{\boldsymbol{\xi}}
\def\bfeta{\boldsymbol{\eta}} \def\bfthz{\boldsymbol{\theta}}
\def\bfdz{\boldsymbol{\delta}}  
\def\bfa{{\bf a}}  \def\bfb{{\bf b}}

\def\bfv{{\boldsymbol{v}}}

\def\bbX{{\mathbb X}}
\def\bbF{{\mathbb F}}
\def\bbH{{\mathbb H}}
\def\bbZ{{\mathbb Z}}
\def\bbN{{\mathbb N}}
\def\bbL{{\mathbb L}}

\def\bbQ{{\mathbb Q}}
\def\bbC{{\mathbb C}}
\def\bbP{{\mathbb P}}
\def\bbO{{\mathbb O}}
\def\bbT{{\mathbb T}}

\def\cP{{\cal P}}
\def\cR{{\cal R}}
\def\cI{{\cal I}}
\def\cS{{\cal S}}
\def\cX{{\cal X}}

\def\cF{{\cal F}}

\def\sF{{\scr F}}
\def\sT{{\scr T}}

\def\sL{{\scr L}}
\def\fn{{\frak n}}
\def\fg{{\frak g}}
\def\cM{{\cal M}}

\def\tbbQ{{\widetilde{\bbQ}}}

\def\lra{\longrightarrow}

\def\lmto{\longmapsto}
\def\ra{\rightarrow}
\def\rmx{{\rm x}}

\def\wt{\widetilde}
\def\calH{{\cal H}}
\def\bfcalH{\boldsymbol{\cal H}}
\def\calK{{\cal K}}
\def\wtH{{\widetilde H}}
\def\qq{{\boldsymbol q}}
\def\frakQ{{\frak Q}}

\def\bfcalC{\boldsymbol{\cal C}}

\def\bfcalK{\boldsymbol{\cal K}}

\begin{document}

\title[Hall polynomials for tame type]{Hall polynomials for tame type}
\author{Bangming Deng and Shiquan Ruan}
\address{Yau Mathematical Sciences Center, Tsinghua University,
Beijing 100084,  China.} \email{bmdeng@@math.tsinghua.edu.cn,\;\;sqruan@@math.tsinghua.edu.cn }
%\address{Department of Mathematical Sciences, Tsinghua University,
%Beijing 100084,  China.} \email{jxiao@@math.tsinghua.edu.cn}

\thanks{Supported partially by the Natural Science Foundation of China.}

%\date{\today}

\subjclass[2000]{17B37, 16G20}
\begin{abstract}
In the present paper we prove that Hall polynomial exists for each triple of decomposition
sequences which parameterize isomorphism classes of coherent sheaves of a domestic weighted projective
line $\bbX$ over finite fields. These polynomials are then used to define the generic Ringel--Hall
algebra of $\bbX$ as well as its Drinfeld double. Combining
this construction with a result of Cramer, we show that Hall polynomials exist for tame quivers,
which not only refines a result of Hubery, but also confirms a conjecture of Berenstein and Greenstein.
\end{abstract}

\maketitle
\begin{center}
{\it To the memory of Professor J.~A. Green}
\end{center}

\bigskip

\section{Introduction}

Inspired by the work of Steinitz \cite{St} and Hall \cite{Hall}, Ringel \cite{Ringel0,Ringel1} introduced
the Hall algebra $H(A)$ of a finite dimensional algebra $A$, whose structure constants are given by the
so-called Hall numbers, and proved that
if $A$ is hereditary and representation finite, then $H(A)$ is isomorphic to the positive part of
the corresponding quantized enveloping algebra. By introducing a bialgebra structure on $H(A)$,
Green \cite{Green} then generalized Ringel's work to arbitrary finite dimensional hereditary algebra $A$
and showed that the composition subalgebra of $H(A)$ generalized by simple $A$-modules
gives a realization of the positive part of the quantized enveloping algebra associated with $A$.
The proof of the compatibility of multiplication and comultiplication on $H(A)$ is based on a marvelous
formula arising from the homological properties of $A$-modules, called Green's formula. We remark
that Lusztig \cite{L91} has obtained a geometric
construction of quantized enveloping algebras in terms of perverse sheaves on representation varieties of quivers.

In case $A$ is representation finite and hereditary, Ringel \cite{Ringel0} showed that the structure constants
of $H(A)$ are actually integer polynomials in the cardinalities of finite fields. The proof is based on
a basic property of the module category of $A$, namely, the directedness of its Auslander--Reiten quiver.
These polynomials are called Hall polynomials as in the classical case;  see \cite{Macdonald}.
Then one can define the generic Hall algebra $H_\qq(A)$ over the polynomial ring
$\bbQ[\qq]$ and its degeneration $H_1(A)$ at $\qq=1$. It was shown by Ringel \cite{R90c} that $H_1(A)$ is isomorphic
to the positive part of the universal enveloping algebra of the semisimple Lie algebra associated with $A$.
Since then, much subsequent work was devoted to the study of Hall polynomials for various classes
of algebras. Recently, Hubery \cite{Hub} provided an elegant proof of the existence of Hall polynomials
for all Dynkin and cyclic quivers by an inductive argument based on Green's formula mentioned above.
Moreover, he proved that Hall polynomials exist for all tame (affine) quivers with respect
to the decomposition classes of Bongartz and Dudek \cite{BD}.

Inspired by the work of Hubery \cite{Hub}, the main purpose of the present paper is to study
Hall polynomials for coherent sheaves of a domestic weighted projective line $\bbX$ over finite fields.
The key idea is again the use of Green's formula. More precisely, we extend the notion of decomposition
classes to that of decomposition sequences, which parameterize isoclasses (isomorphism classes) of coherent sheaves of
$\bbX$ over finite fields, and show that Hall polynomial exists for each triple of decomposition
sequences. These polynomials are then applied to define an algebra $H_\bfv(\bbX)$ which is a free
module over the Laurent polynomial ring $\bbQ[\bfv,\bfv^{-1}]$ with a basis all the
decomposition sequences. By extending $H_\bfv(\bbX)$ via formally adding certain elements
constructed in \cite{BurSch}, we obtain the generic Ringel--Hall algebra $\bfcalH_\bfv(\bbX)$ of $\bbX$.
By further introducing Green's pairing on $\bfcalH_\bfv(\bbX)$,
we construct its Drinfeld double $D\bfcalH_\bfv(\bbX)$ over $\bbQ(\bfv)$. Combining
this construction with \cite[Prop.~5]{cram}, we show that Hall polynomials exist for decomposition sequences
associated with a tame quiver. This result refines the main theorem of Hubery \cite{Hub} and also
confirms a conjecture of Berenstein and Greenstein \cite[Conj.~3.4]{BG}.

The paper is organized as follows. Section 2 gives a brief introduction on the category of coherent
sheaves over a weighted projective line $\bbX$ and recalls the definition of the Ringel--Hall algebra of
$\bbX$ over a finite field and the Green's formula as well. In Section 3 we define decomposition sequences for a domestic
weighted projective line and give some preparatory results which are needed in Section 4 to prove the
existence of Hall polynomials. Section 5 is devoted to defining the generic Hall algebra
$\bfcalH_v(\bbX)$ of $\bbX$ as well as its Drinfeld double $D\bfcalH_v(\bbX)$. In the final section, we show that Hall
polynomials exist for tame quivers.

\section{The category of coherent sheaves over a weighted projective line}

In this section we review the category of coherent sheaves over a
weighted projective line and its basic properties, and we also introduce Hall algebras
and Green's formula. For further fundamental concepts and facts on categories of coherent
sheaves over weighted projective lines and on Hall algebras,
we refer to \cite{GeigleLenzing,CH} and \cite{Ringel4,Sch2006,DDPW}.

\subsection{The category of coherent sheaves}

Let $k$ be an arbitrary field. A \emph{weighted projective line}
$\bbX=\bbX_k$ over $k$ is specified by giving a \emph{weight
sequence} ${\bf p}=(p_{1},p_{2},\ldots, p_{t})$ of positive
integers, and a collection
${\boldsymbol\lambda}=(\lambda_{1},\lambda_{2},\ldots, \lambda_{t}) $ of
distinct points in the projective line $\bbP^{1}(k)$ which can be
normalized as $\lambda_{1}=\infty, \lambda_{2}=0, \lambda_{3}=1$.
More precisely, let $\bbL=\bbL(\bf p)$ be the rank one abelian group with
generators $\vec{x}_{1}, \vec{x}_{2}, \ldots, \vec{x}_{t}$ and the
relations
\[ p_{1}\vx_1=p_{2}\vx_2=\cdots=p_{t}\vx_t=:\vec{c},\]
where $\vec{c}$ is called the \emph{canonical element} of
$\mathbb{L}$. Denote by $S$ the commutative algebra
\[S=k[X_{1},X_{2},\cdots,
X_{t}]/{\frak a}:= k[{\rm x}_{1},{\rm x}_{2}, \ldots, {\rm x}_{t}],\] where
${\frak a}=(f_{3},\ldots,f_{t})$ is the ideal generated by
$f_{i}=X_{i}^{p_{i}}-X_{2}^{p_{2}}+\lambda_{i}X_{1}^{p_{1}},
i=3,\ldots, t$. Put $I=\{1,2,\ldots,t\}$. Then $S$ is $\mathbb{L}$-graded by setting
$$\mbox{deg}(\rmx_{i})=\vx_i\; \text{ for each $i\in I$.}$$
Moreover, each element $\vx\in\bbL$ has the normal form $\vx=\sum_{i\in I} l_i\vx_i+l\vc$ with
$0\leq l_i<p_i$ and $l\in\bbZ$. We denote by $\bbL_+$ the positive cone of $\bbL$ which consists of
those $\vx$ with $l\geq 0$. Finally, the weighted projective line associated with $\bf p$ and $\boldsymbol\lz$
is defined to be $\bbX=\rm{Spec}^{\bbL}S$.

According to \cite{GeigleLenzing}, the set of nonzero prime homogeneous
elements in $S$ is partitioned into two sets: the
exceptional primes $\rmx_1,\ldots, \rmx_t$ and the ordinary primes
$f(\rmx_{1}^{p_{1}}, \rmx_{2}^{p_{2}})$, where $f(T_1, T_{2})$ is a prime
homogeneous polynomial in $k[T_1,T_2]$ which are distinct from $T_1, T_2$ and
$T_2-\lz_iT_1$ for $i\in I$. The exceptional primes
correspond to the points $\lz_1, \ldots, \lz_t$, called
exceptional points and denoted by $x_1,\ldots,x_t$, respectively, while the ordinary primes correspond to the
remaining closed points of $\bbP^{1}(k)$, called ordinary
points. For convenience, we denote by ${\mathbb H}_k$ the set of ordinary
points. For each $z\in {\mathbb H}_k$, its degree $\deg(z)$ is defined to be the degree
of the corresponding prime homogeneous polynomial.

The category of coherent sheaves on $\bbX$ can be defined as the
quotient of the category of finitely generated $\mathbb{L}$-graded
$S$-modules over the Serre subcategory of finite length modules, that is,
$$\coh\bbX:={\rm mod}^{\mathbb{L}}(S)/\mbox{mod}_{0}^{\mathbb{L}}(S).$$
The free module $S$ gives the structure sheaf $\co$. Each line
bundle is given by the grading shift $\co(\vec{x})$ for a uniquely
determined element $\vec{x}\in \mathbb{L}$, and there is an
isomorphism
\[\Hom(\co(\vec{x}), \co(\vec{y}))\cong S_{\vec{y}-\vec{x}}.\]
Moreover, $\coh\bbX$ is a hereditary abelian category with Serre
duality of the form
\[D\Ext^1(X, Y)\cong\Hom(Y, X(\vec{\omega})),\]
where $D=\Hom_k(-,k)$, and $\vec{\omega}:=(t-2)\vec{c}-\sum_{i\in I}\vec{x}_{i}$ is
called the \emph{dualizing element} of $\mathbb{L}$. This implies
the existence of almost split sequences in $\coh\bbX$ with the
Auslander--Reiten translation $\tau$ given by the grading shift with
$\vec{\omega}$.

Recall that $\coh\bbX$ admits a splitting torsion pair $({\rm coh}_0\mbox{-}\bbX,\vect{\bbX})$,
where ${\rm coh}_0\mbox{-}\bbX$ and $\vect{\bbX}$ are full subcategories
of torsion sheaves and vector bundles, respectively. Moreover, ${\rm coh}_0\mbox{-}\bbX$ decomposes as a
direct product of orthogonal tubes
$${\rm coh}_0\mbox{-}\bbX=\prod_{z\in\bbH_k}{\rm
coh}_z\mbox{-}\bbX\times\prod_{i\in I}{\rm coh}_i\mbox{-}\bbX$$
 where each ${\rm coh}_z\mbox{-}\bbX$ is a homogeneous tube, which is equivalent to the category of
nilpotent representations of the Jordan quiver over the residue
field $k_z$, while each ${\rm coh}_i\mbox{-}\bbX$ is a non-homogeneous tube, which is equivalent to the
category of nilpotent representations of the cyclic quiver with $p_i$ vertices. By a classical result,
the isoclasses (isomorphism classes) of objects in ${\rm coh}_z\mbox{-}\bbX$ are indexed by partitions,
while those in ${\rm coh}_i\mbox{-}\bbX$
are indexed by multipartitions; see, for example, \cite{Ringel3}. More precisely, for each $z\in\bbH_k$,
${\rm coh}_z\mbox{-}\bbX$ admits a unique simple object $S_z$ and, up to isomorphism, each object in
${\rm coh}_z\mbox{-}\bbX$ has the form $S_k(\pi,z)=\oplus_{r=1}^sS_z[\pi_r]$, where
$\pi=(\pi_1,\ldots,\pi_s)$ is a partition and $S_z[\pi_r]$ is the unique indecomposable object of length
$\pi_r$. While for each $i\in I$, there are $p_i$ simple objects $S_{i,0},\ldots, S_{i,p_i-1}$
in ${\rm coh}_i\mbox{-}\bbX$. For each $0\leq j\leq p_i-1$ and $l\geq 1$, let $S_{i,j}[l]$ denote
the indecomposable object in ${\rm coh}_i\mbox{-}\bbX$ of length $l$ with top $S_{i,j}$.

It is known that the Grothendieck group $K_0(\bbX)$ of $\coh\bbX$ is a free abelian group
with a basis $\co(\vec{x})$ with $0\leq\vec{x}\leq \vec{c}$, where we still
write $X\in K_0(\bbX)$ for the isoclass of an object $X\in \coh\bbX$.
Let $p={\rm l.c.m.}(p_{1},\ldots, p_{t})$ be the least common multiple of $p_1,\ldots,p_t$ and
$\delta: \mathbb{L}\to \mathbb{Z}$ be the homomorphism defined by
$\delta(\vec{x}_{i})=\frac{p}{p_{i}}.$ The \emph{determinant} map is the group homomorphism
$$\det: K_0(\bbX)\lra {\mathbb L},\;{\co}(\vec{x})\longmapsto \vec{x}.$$
The rank function on $K_0(\bbX)$ is given by the rule
$\mbox{rk}(\co(\vec{x}))=1$ while the degree function is given by
the rule $\deg(\co(\vec{x}))=\dz(\vx)$. For each non-zero object
$X\in \coh\bbX$, define the \emph{slope} of $X$ as $\mu(X)=\frac{\mbox{deg}(X)}{\mbox{rk}(X)}.$
The Euler form on $K_0(\bbX)$ is given by
\[\langle X,Y\rangle=\dim_{k}\Hom(X,Y)-\dim_{k}\Ext^{1}(X,Y).\]
for any $X,Y\in\coh\bbX$. Its symmetrization is defined by
$$(X,Y)=\langle X,Y\rangle +\langle Y,X\rangle.$$

In the present paper we mainly focus on weighted projective lines $\bbX$ of domestic type,
i.e., $\dz(\vw)<0$. In this case, the Auslander--Reiten quiver $\Gamma(\vect\bbX)$ of $\vect\bbX$ consists of a single standard
component of the form $\bbZ\widetilde \Delta$, where $\widetilde\Delta$ is an extended Dynkin
diagram associated with the weight sequence $\bf p$. Moreover, the full subcategory of indecomposable
vector bundles on $\bbX$ is equivalent to the mesh category of  $\Gamma(\vect\bbX)$. Furthermore,
for any two indecomposable objects $X, Y\in\coh\bbX$, $\mbox{Hom}(X,Y)\neq 0$ implies $\mu(X)\leq \mu(Y).$

The following result will be needed later on.

\begin{Lem}\label{vector bundle factorization} Let $E$ be an indecomposable vector bundle. Then
there is an exact\linebreak sequence $0\to L\to E\to F\to 0$ in $\vect\bbX$ such that $L$ is
a line bundle and $\Ext^1(F, L)\cong k$.
\end{Lem}

\begin{pf} Choose a line bundle $L$ of maximal degree such that $\Hom(L, E)\neq 0$. Then we get an exact sequence
\begin{equation}\label{fitting of E}
0\lra L\lra E\lra F\lra 0
 \end{equation}
 in $\coh\bbX$. We claim that $F$ is a vector bundle. Otherwise, there exists a simple subsheaf $S$ of
 $F$, which yields the following pullback commutative diagram:
$$ \xymatrix{
  0\ar[r] & L \ar@{=}[d] \ar[r] &  \ar[d] \ar[r]L' & S \ar[d]\ar[r]&0 \\
  0\ar[r] & L\ar[r] & E \ar[r] & F\ar[r]&0.  }$$
This gives a line bundle $L'$ satisfying that $\Hom(L', E)\neq 0$ and $\deg L'>\deg L$, contradicting the choice of $L$.

Applying $\Hom(-,L)$ to the exact sequence (\ref{fitting of E}) gives the exact sequence
$$0\lra \Hom(L,L)\lra \Ext^1(F,L)\lra \Ext^1(E,L).$$
Note that $\Hom(L,L)\cong k$ and $\Ext^1(E,L)\cong D\Hom(L(-\vw),
E)=0$ since $\deg L(-\vw)>\deg L$. Therefore, $\Ext^1(F,L)\cong k$.
\end{pf}

\subsection{The Hall algebra of coherent sheaves}

Let $k$ be a finite field. Given objects $Z, X_1,\ldots, X_t$ in $\coh\bbX_k$, define
$F^Z_{X_1,\ldots, X_t}$ to be the number of filtrations
$$Z=Z_0\supseteq Z_1\supseteq \cdots \supseteq Z_{t-1}\supseteq Z_t=0$$
such that $Z_{s-1}/Z_s\cong X_s$ for all $1\le s\le t$, called the
{\it Hall number} associated with $Z,X_1,\ldots,X_t$.

For each object $X\in\coh\bbX_k$, put $a_X=|\Aut(X)|$, the cardinality of the automorphism
group $\Aut(X)$ of $X$. The following result is taken from \cite{Riedtmann,Peng}.

\begin{Lem} \label{Riedtmann-Peng-formula}
Let $X, Y, Z$ be three objects in $\coh\bbX_k$. Then
$$F^Z_{X,Y}=\frac{|\Ext^1(X,Y)_Z|}{|\Hom(X,Y)|}\cdot \frac{a_Z}{a_Xa_Y},$$
 where $\Ext^1(X,Y)_Z$ denotes the subset of $\Ext^1(X,Y)$ consisting of equivalence classes
 of exact sequences in $\coh\bbX$ of the form $0\ra Y\ra Z\ra X\ra 0.$
\end{Lem}

Now let $k$ be a finite field with $q$ elements and let $v_q$ denote the square root
$\sqrt{q}$ of $q$. For each $M\in\coh\bbX_k$, let $[M]$ denote the isoclass of $M$.
By definition, the Ringel--Hall algebra $H(\bbX_k)$ of the category of coherent
sheaves on $\bbX_k$ is the free module over the ring $\bbQ[v_q,v_q^{-1}]$
with basis $\{[M]\mid M\in\coh\bbX_k\}$, and the multiplication is given by
$$[M][N]=v_q^{\langle M, N\rangle}\sum\limits_{[R],\,R\in\text{coh-}\bbX_k}F_{M ,N}^{R} [R].$$
 By a result of Green \cite{Green}, $H(\bbX_k)$ is a bialgebra with comultiplication defined by
$$\Dz_k([R])=\sum\limits_{[M],[N]}v_q^{\langle M,
N\rangle}F_{M ,N}^{R}\frac{a_{M}a_{N}}{a_{R}}[M]\otimes [N].$$
 In fact, the associativity of multiplication and the
coassociativity of comultiplication follow from the identity
$$\sum\limits_{X}F_{A ,B}^{X}F_{X ,C}^{S}=\sum\limits_{X}F_{A ,X}^{S}F_{B ,C}^{X},$$
 where the sums on both sides are actually taken over isoclasses of
objects in $\coh\bbX_k$, though we shall often use this more
convenient notation. Furthermore, the compatibility of multiplication and
comultiplication is encoded in following marvellous
formula---the so-called Green's formula---which plays a fundamental
role in the study of Hall algebras.

\begin{Lem} [\cite{Green}]  \label{Green-formula-Lem} For each quadruple $(M,N,X,Y)$ of objects in $\coh\bbX_k$,
we have the equality
\begin{equation}\label{Green-formula}
\sum_{E}F_{M,N}^EF_{X,Y}^E/a_E = \sum_{A,B,C,D}q^{-\langle A,D\rangle}F_{A,B}^MF_{C,D}^NF_{A,C}^XF_{B,D}^Y
\frac{a_Aa_Ba_Ca_D}{a_Ma_Na_Xa_Y}.
\end{equation}
\end{Lem}

Following an idea in \cite[Sect.~3]{Hub}, if $\Ext^1(X,Y)=0$, then the left-hand side of
\eqref{Green-formula} contains only one term
$$F_{M,N}^{X\oplus Y} F_{X,Y}^{X\oplus Y}/a_{X\oplus Y},$$
 which, by Lemma \ref{Riedtmann-Peng-formula}, is equal to
$$F_{M,N}^{X\oplus Y}\frac{a_{X\oplus Y}}{|\Hom(X,Y)| a_Xa_Y}\cdot\frac{1}{a_{X\oplus Y}}
=q^{-\langle X,Y\rangle}\frac{1}{a_Xa_Y}F_{M,N}^{X\oplus Y}.$$
 Thus, in this case, Green's formula \eqref{Green-formula} is simplified to the form
\begin{equation}\label{simple-Green-formula}
F_{M,N}^{X\oplus Y}=\sum_{A,B,C,D}q^{\langle X,Y\rangle-\langle A,D\rangle}
F_{A,B}^MF_{C,D}^NF_{A,C}^{X}F_{B,D}^{Y}\frac{a_Aa_Ba_Ca_D}{a_Ma_N}.
\end{equation}

\subsection{The elements $\Thz_{\vx}$, $T_r$ and $Z_r$ in $H(\bbX_k)$}

As above, let $k$ be a finite field with $q$ elements. In the following we recall from \cite{BurSch}
the definition of some special elements in the Ringel--Hall algebra $H(\bbX_k)$ which will be needed later on.

By \cite[5.5]{BurSch}, for each $\vx\in\bbL_+$ with normal form $\vx=\sum_{i\in I} l_i\vx_i+l\vc$,
define the element $\Thz_{\vx}=\Thz_{\vx,q}$ via the formula
\begin{equation}\label{def of theta}\Dz([\co])=[\co]\otimes 1+\sum\limits_{\vx\in\bbL_+}\Thz_{\vx}\otimes [\co(-\vx)].
\end{equation}
Then $\Thz_{\vx}$ can be written as
\begin{equation}\label{def-Theta_x}
\aligned
{}\qquad\Thz_{\vx}&=v_q^{l+m }\sum\limits_{z_j,n_j,m_i}\prod\limits_{j}(1-v_q^{-2\deg(z_j)})\times\\
&\qquad \prod_{i\in I,(m_i,l_i)\neq (0,0)}(1-v_q^{-2})
\big[\bigoplus_{j}S_{z_j}[n_j]\oplus\bigoplus_{i\in I}S_{i,0}[m_ip_i+l_i]\big],
\endaligned
\end{equation}
where $m=| \{i\mid l_i\neq 0\}|$, and the sum ranges over
all tuples of distinct ordinary points $z_j$ and nonnegative integers $n_j$, $m_i$ satisfying
$\sum_{j}n_j\deg(z_j)+\sum_{i\in I}m_i=l$.

However, the definition of the elements $T_r=T_{r,q}$ and $Z_r=Z_{r,q}$
in $H(\bbX_k)$ for all $r\geq 1$ are rather complicated, so we refer to \cite[Sect.~6]{BurSch}.
We emphasize that the $Z_r$ commute with $[S]$ for all torsion sheaves $S$ in $\coh \bbX_k$.

\section{Hall polynomials for a domestic weighted projective line}

In this section, we introduce the notion of Hall polynomials and
prove that Hall polynomials exist for a domestic weighted projective line.

As in the previous section, let $\bbX=\bbX_k$ be a domestic weighted
projective line over a finite field $k$. By $\chi=\chi(\bbX)$ we denote
the set of isoclasses of objects in
$\coh\bbX$ which clearly depends on the ground field $k$. Let
$\chi_t$ and $\chi_f$ be the subsets of $\chi$ consisting of the
isoclasses of torsion sheaves and vector bundles, respectively.
Further, let $\chi_\nh$ be the subset formed by the isoclasses of
sheaves without homogeneous regular summands. In other words,
$\chi_\nh$ consists of isoclasses of sheaves whose indecomposable
summands are either vector bundles or torsion sheaves lying in
non-homogeneous tubes. Hence, the set $\chi_\nh$ can be described
combinatorially and is independent of $k$. Moreover, each sheaf in a
homogeneous tube corresponding to a point in $\bbH_k$ is determined
by a partition. For each $\az\in\chi$, we fix a representative
$S_k(\az)$ in the class $\az$. Given $\az,\bz\in\chi$, we write
$\az\oplus\bz$ for the isoclass of $S_k(\az)\oplus S_k(\bz)$. Thus,
each $\az\in\chi$ can be uniquely decomposed as
$\az=\az_t\oplus\az_f$ with $S_k(\az_t)\in {\rm coh}_0\mbox{-}\bbX$
and $S_k(\az_f)\in\vect\bbX$. % Note that $\rk S_k(\az)$ is
%independent of the field $k$. We set $\rk\gz=\rk S_k(\gz,\ul{x})$.

\subsection{Segre sequences and Hall polynomials}
A {\it Segre sequence} is a sequence $\lz=((\lz^{(1)},d_1), (\lz^{(2)},d_2), \ldots, (\lz^{(r)},d_r))$ of
pairs $(\lz^{(i)},d_i)$, where $\lz^{(i)}$ are partitions and $d_i$ are positive integers with
$d_1\leq d_2\leq \cdots\leq d_r$. Such a sequence is said to be of type $\ul{d}=(d_1, d_2, \ldots, d_r)$.
If all $\lz^{(i)}$ are the empty partition, then we simply write $\lz=\emptyset$. A {\it decomposition
sequence of type $\ul{d}$} is by definition a pair $\bfaz=(\alpha, \lz)$, where $\alpha\in\chi_\nh$
and $\lz$ is a Segre sequence of type $\ul{d}$.

\begin{Rem} \label{equiv-Segre-class} For a partition $\pi$ and a positive integer $d$, by inserting the pair
$(\pi,d)$ to a Segre sequence $\lz=((\lz^{(1)},d_1), (\lz^{(2)},d_2), \ldots, (\lz^{(r)},d_r))$ we mean
the Segre sequence
$$\mu:=((\lz^{(1)},d_1),\ldots, (\lz^{(i)},d_i), (\pi,d), (\lz^{(i+1)},d_{i+1}) \ldots, (\lz^{(r)},d_r)),$$
where $1\leq i\leq r$ satisfies $d_i\leq d<d_{i+1}$. In this case, we also say that $\lz$ is obtained from
$\mu$ by removing the pair $(\pi,d)$. In particular, any finitely many Segre sequences can be
converted to Segre sequences of same type via inserting and removing some pairs $(\emptyset, d)$.
Finally, two decompositions sequences $\bfaz=(\alpha, \lz)$ and $\bfbz=(\bz, \mu)$ are identified if
$\az=\bz$ and $\lz$ can be obtained from $\mu$ by inserting and removing some pairs $(\emptyset, d)$.
\end{Rem}

For a finite field $k$, denote by ${\cal X}_k(\ul{d})$ the set of
sequences $\ul{z}=(z_1, \ldots, z_r)$ of pairwise distinct points in
$\bbH_k$ with $\deg(z_i)=d_i$. Note that the cardinality of ${\cal X}_k(\ul{d})$ depends on the ground
field $k$, and ${\cal X}_k(\ul{d})$ is possibly empty. However, ${\cal X}_k(\ul{d})\not=\emptyset$ when $|k|\gg0$.

For each Segre sequence $\lz$ of type $\ul{d}$ and $\ul{z}=(z_1, \ldots, z_r)\in {\cal X}_k(\ul{d})$, define
$$S_k(\lz, \ul{z}):=\bigoplus\limits_{i=1}^{r}S_k(\lz^{(i)}, z_i)\in\coh\bbX_k,$$
where $S_k(\lz^{(i)}, z_i)\in{\rm coh}_{z_i}\text{-}\bbX_k$ is determined by the partition $\lz^{(i)}$.
%called a torsion sheaf of type $\lz$.
Further, for each decomposition sequence $\bfaz=(\alpha, \lz)$ of
type $\ul{d}$, define
$$S_k(\bfaz, \ul{z}):=S_k(\alpha)\oplus S_k(\lz, \ul{z}).$$
If $\lz=\emptyset$, then $\bfaz=(\az,\emptyset)\in\chi_\nh$ and
$S_k(\bfaz)=S_k(\az)$. Clearly, $S_k(\bfaz, \ul{z})\in\vect\bbX_k$ if
and only if $\az=\az_f\in\chi_f$ and $\lz=\emptyset$. In this case,
$\bfaz$ is said to be of {\it torsion-free type} (we also simply
write $\bfaz\in\chi_f$), and we write
$S_k(\bfaz, \ul{z})=S_k(\bfaz)=S_k(\az)$. Furthermore, $S_k(\bfaz,
\ul{z})\in\tor\bbX_k$ if and only if $\az\in \chi_\nh\cap \chi_t$. In
this case, $\bfaz$ is said to be of {\it torsion type}.

Given $\bfaz=(\alpha, \lz)$ of type $\ul{d}=(d_1,\ldots,d_r)$, if $\mu$ is obtained from $\lz$ by
removing a pair $(\emptyset,d_s)$, where $1\leq s\leq r$
and $d_s<d_{s+1}$, then $\mu$ is of type $\ul{d'}=(d_1,\ldots,d_{s-1},d_{s+1},\ldots,d_r)$ and, moreover,
for each $\ul{z}=(z_1, \ldots, z_r)\in {\cal X}_k(\ul{d})$,
$$S_k(\bfaz, \ul{z})\cong S_k((\az,\mu), \ul{z'}),$$
 where $\ul{z'}=(z_1,\ldots,z_{s-1},z_{s+1},\ldots,z_r)\in {\cal X}_k(\ul{d'})$. Therefore,
by setting $\bfbz=(\az,\mu)$, the two sets
$$\{S_k(\bfaz, \ul{z})\mid \ul{z}\in {\cal X}_k(\ul{d})\}\;\text{ and }\; \{S_k(\bfbz, \ul{z'})\mid \ul{z'}\in {\cal X}_k(\ul{d'})\}$$
 give rise to the same family of isoclasses of coherent sheaves in $\coh\bbX_k$.

By $\cS$ we denote the set of all decomposition sequences (up to the identification
in Remark \ref{equiv-Segre-class}) and by
$\cS_t$ (resp., $\cS_f$) its subset of decomposition sequences of
torsion type (resp., torsion-free type). Note that $\cS_f$ can be
identified with $\chi_f$ which is a subset of $\chi_\nh$. Now for each
$\bfaz=(\az,\lz)\in\cS$, the decomposition $\az=\az_t\oplus\az_f$
gives two decomposition sequences
$$\bfaz_t=(\az_t,\lz)\in\cS_t \;\text{ and }\;\bfaz_f=(\az_f,\emptyset)\in\cS_f$$
such that for each $\ul{z}\in\cX_k(\ul{d})$,
$$S_k(\bfaz_t,\ul{z})=S_k(\bfaz, \ul{z})_t\;\text{ and }\; S_k(\bfaz_f,\ul{z})=S_k(\bfaz, \ul{z})_f=S_k(\az_f).$$
We simply write $\bfaz=\bfaz_t\oplus\bfaz_f$.

Give $\bfaz,\bfbz\in\cS$ of type $\ul{d}$, it is easy to see from the definition that the value
$\langle S_k(\bfaz,\ul{z}), S_k(\bfbz,\ul{z})\rangle$ is a constant for any field $k$ with
$|k|\gg0$ and $\ul{z}\in{\cal X}_k(\ul{d})$. Thus, we simply
put
$$\langle\bfaz, \bfbz\rangle=\langle S_k(\bfaz,\ul{z}), S_k(\bfbz,\ul{z})\rangle\;\text{ and }\;
(\bfaz, \bfbz)=\langle\bfaz, \bfbz\rangle+\langle\bfbz, \bfaz\rangle.$$

\begin{Defn} \label{Def-Hall-poly} Given $\bfaz, \bfbz,\bfgz\in\cS$ of type $\ul{d}$,
if there exists a polynomial $\vphi_{\bfaz,\bfbz}^{\bfgz}\in\bbZ[T]$
such that for each finite field $k$ of $q$ elements with $q\gg0$,
$$\vphi_{\bfaz,\bfbz}^{\bfgz}(q)= F_{S_k(\bfaz, \ul{z}), S_k(\bfbz, \ul{z})}^{S_k(\bfgz, \ul{z})}
\quad\textrm{for all }\ul{z}\in \cX_k(\ul{d}),$$ then we say that
the Hall polynomial $\vphi_{\bfaz,\bfbz}^{\bfgz}$ exists for
$\bfaz$, $\bfbz$ and $\bfgz$.
\end{Defn}

One of our main purposes in the present paper is to prove the following result.

\begin{Thm} \label{existence-HallPoly} For arbitrary $\bfaz, \bfbz,\bfgz\in\cS$ of type $\ul{d}$,
the Hall polynomial $\vphi_{\bfaz,\bfbz}^{\bfgz}$ exists.
\end{Thm}

The theorem will be proved in the next section. In the following we present some preparatory results
which will be needed for the proof of the main theorem.

\begin{Lem}[\cite{Ringel0}]\label{Ringel dividing lemma}
Let $\phi,\psi\in\bbZ[T]$ and assume $\psi$ is monic. Then $\psi$ divides $\phi$ if
and only if the integer $\psi(q)$ divides the integer $\phi(q)$ for infinitely many $q\in\mathbb{Z}$.
\end{Lem}

\begin{Lem}\label{poly for hom}
Given decomposition sequences $\bfaz$ and $\bfbz$ of type $\ul{d}$,
there exists a monic integer polynomial $h_{\bfaz,\bfbz}\in\bbZ[T]$
such that for any finite field $k$ of $q$ elements with $q\gg0$,
$$h_{\bfaz,\bfbz}(q)=|\Hom(S_k(\bfaz, \ul{z}), S_k(\bfbz, \ul{z}))| \quad\textrm{for all } \ul{z}\in \cX_k(\ul{d}).$$
\end{Lem}

\begin{pf}
Assume $\bfaz=(\az, \lz)$ and $\bfbz=(\bz, \mu)$ with $\az,
\bz\in\chi_\nh$ and $\lz, \mu$ are Segre sequences. Since the
$\Hom$-functor commutes with direct sum, it suffices to prove the
lemma when both $S_k(\bfaz, \ul{z})$ and $S_k(\bfbz, \ul{z})$ are
indecomposable. Thus, $\az=0$ or $\lz=\emptyset$, and $\beta=0$ or $\mu=\emptyset$.
We treat each of the four cases as follows:

\begin{itemize}
\item[(i)] {\it Case $\lz=\emptyset$ and $\mu=\emptyset$.} Then the assertion
follows from the fact that $\chi_\nh$ can be described
combinatorially.

\item[(ii)] {\it Case $\az=0$ and $\bz=0$.} Then $S_k(\bfaz, \ul{z})$ and
$S_k(\bfbz,\ul{z})$ are $\Hom$-orthogonal or belong to the same tube. This case
follows from the representation theory of a cyclic quiver; see \cite{Ringel3,Guo}.

\item[(iii)] {\it Case $\lz=\emptyset$ and $\bz=0$.} Then $\Ext^1(S_k(\bfaz, \ul{z}), S_k(\bfbz, \ul{z}))=0$.
By Riemann--Roch formula in \cite{GeigleLenzing},
$$\qquad\quad \rk \bfaz \deg \bfbz=\sum\limits_{i=1}^p\langle
S_k(\bfaz, \ul{z}), \tau^i S_k(\bfbz, \ul{z})\rangle = p\dim_k
\Hom(S_k(\bfaz, \ul{z}), S_k(\bfbz, \ul{z})).$$
Thus,
$$\dim_k\Hom(S_k(\bfaz, \ul{z}), S_k(\bfbz, \ul{z}))=\frac{1}{p}\rk\bfaz \deg \bfbz,$$
which is an integer (since $p$ divides $\deg
\bfbz$) and independent of the field $k$, as desired.

\item[(iv)] {\it Case $\az=0$ and $\mu=\emptyset$}. This case can be proved by
an argument similar to that in case (iii).
\end{itemize}
\end{pf}

The following is an easy consequence of the lemma above.

\begin{Lem}\label{poly for auto}
Given a decomposition sequence $\bfaz$ of type $\ul{d}$, there exists
a monic integer polynomial $a_{\bfaz}\in\bbZ[T]$ such that, for any
finite field $k$ of $q$ elements with $q\gg0$,
$$a_{\bfaz}(q)=|\Aut(S_k(\bfaz, \ul{z}))| \quad\textrm{for all }\ul{z}\in \cX_k(\ul{d}).$$
\end{Lem}

\begin{Prop}\label{extistence-vb} If $\bfaz,\bfbz,\bfgz\in\cS_f$, then the Hall polynomial
$\vphi_{\bfaz,\bfbz}^{\bfgz}$ exists. If, moreover, $\Ext^1(S_k(\bfaz),S_k(\bfbz))\cong k$, then
$\vphi_{\bfaz,\bfbz}^{\bfgz}$ is monic.
\end{Prop}

\begin{pf} By the assumption, we can write $\bfaz=(\az,\emptyset)$, $\bfbz=(\bz,\emptyset)$,
and $\bfgz=(\gz,\emptyset)$ for $\az,\bz,\gz\in\chi_f$.

Let $k$ be a finite field and take a complete slice $\scr{S}$
in $\vect\bbX_k$, which gives a tilting bundle $T$, such that every
indecomposable direct summand of $S_k(\az), S_k(\bz)$ and $S_k(\gz)$
is generated by $T$. For instance, the slice $\scr{S}$ can be
taken such that each indecomposable direct summand of $S_k(\az),
S_k(\bz)$ and $S_k(\gz)$ lies on the right hand side of the slice in
the Auslander--Reiten quiver of $\vect\bbX_k$. It is well known that
the endomorphism algebra $\Lambda:=\End(T)$ is tame hereditary and
$\cF:=\Hom(T, -)$ induces an equivalence between certain exact full
subcategories of $\coh\bbX_k$ and $\mbox{mod-}\Lambda$. In
particular, the images of $S_k(\az), S_k(\bz)$ and $S_k(\gz)$ under
$\cF$ belong to the preprojective component of $\mbox{mod-}\Lambda$.
Then the existence of $\vphi_{\az,\bz}^{\gz}$ follows from the fact
that Hall polynomials exists for any three preprojective
$\Lambda$-modules. (Note that the latter can be proved by using the
arguments similar to those in \cite[Th.~1]{Ringel0}.)

Now assume $\Ext^1(S_k(\az),S_k(\bz))\cong k$. Thus,
$$|\Ext^1(S_k(\az),S_k(\bz))_{S_k(\gz)}|=0,\;1,\;\text{or}\; |k|-1.$$
By Lemma \ref{Riedtmann-Peng-formula}, we have
$$F_{S_k(\az), S_k(\bz)}^{S_k(\gz)}=\frac{|\Ext^1(S_k(\az), S_k(\bz))_{S_k(\gz)}|}{|\Hom(S_k(\az), S_k(\bz))|}
\frac{|\Aut(S_k(\gz))|}{|\Aut(S_k(\az))||\Aut(S_k(\bz))|}.$$
Then each term on the right hand side of the above equality
is a monic integer polynomial in $|k|$. This implies that $\vphi_{\az,\bz}^{\gz}$ is monic.
\end{pf}

\begin{Prop} \label{exist-regular} Let $\bfaz=(\az,\lz)$, $\bfbz=(\bz,\mu)$ and $\bfgz=(\gz,\nu)$ be
decomposition sequences of type $\ul{d}=(d_1, \ldots, d_r)$. If
$\bfaz, \bfbz, \bfgz\in\cS_t$, then the Hall polynomial
$\vphi_{\bfaz,\bfbz}^{\bfgz}$ exists. Moreover, for fixed
decomposition sequences $\bfaz$ and $\bfbz$ of torsion type, there
are only finitely many decomposition sequences $\bfgz$, and vice
versa, such that $\vphi_{\bfaz,\bfbz}^{\bfgz}\not=0$.
\end{Prop}

\begin{pf} Let $k$ be a finite field. Since there are no nonzero homomorphisms between objects
in distinct tubes of $\coh\bbX_k$, we have for each $\ul{z}=(z_1,\ldots,z_r)\in\cX_k(\ul{d})$,
$$\aligned
F_{S_k(\bfaz,\ul{z}),S_k(\bfbz,\ul{z})}^{S_k(\bfgz,\ul{z})}&=F_{S_k(\alpha),S_k(\beta)}^{S_k(\gamma)}
F_{S_k(\lz, \ul{z}), S_k(\mu, \ul{z})}^{S_k(\nu, \ul{z})}\\
&=F_{S_k(\alpha),S_k(\beta)}^{S_k(\gamma)}
\prod_{i=1}^{r} F_{S_k(\lz^{(i)}, z_i),S_k(\mu^{(i)}, z_i)}^{S_k(\nu^{(i)}, z_i)}.
\endaligned$$
By the existence of Hall polynomials for nilpotent representations of cyclic quivers (including
the classical case; see \cite{Ringel3,Guo,Hall, Macdonald}),
the Hall polynomials $\vphi_{\alpha, \beta}^{\gamma}(T)$ and $\vphi_{\lz^{(i)},\mu^{(i)}}^{\nu^{(i)}}(T)$
($1\leq i\leq r$) exist. Therefore, the polynomial
$$\vphi_{\alpha, \beta}^{\gamma}(T)\prod_{i=1}^{r}\vphi_{\lz^{(i)},\mu^{(i)}}^{\nu^{(i)}}(T^{d_i})\in\bbZ[T]$$
is the required Hall polynomial $\vphi_{\bfaz,\bfbz}^{\bfgz}(T)$.

The second assertion follows from the properties of Hall polynomials for nilpotent representations of a cyclic quiver.
\end{pf}

\section{Proof of Theorem \ref{existence-HallPoly}}

This section is devoted to proving Theorem \ref{existence-HallPoly}.
In the following we always assume that $\bfaz=(\az,\lz),
\bfbz=(\bz,\mu)$, and $\bfgz=(\gz,\nu)$ are decomposition sequences
of type $\ul{d}$. Further, for a finite field $k$
and $\ul{z}\in\cX_k(\ul{d})$, we will simply put
$$M=M_k(\ul{z}):=S_k(\bfaz, \ul{z}),\; N=N_k(\ul{z}):=S_k(\bfbz, \ul{z}),\;\text{ and }\; Z=Z_k(\ul{z}):=S_k(\bfgz, \ul{z}).$$

\subsection{Reduction 1} To prove Theorem
\ref{existence-HallPoly}, it suffices to prove the existence of Hall
polynomials $\vphi_{\bfaz,\bfbz}^{\bfgz}$ for all $\bfaz,\bfbz,
\bfgz\in\cS$ with $\bfgz\in \cS_f=\chi_f$.

\begin{pf} If $\bfgz$ is of torsion type and $F_{M_k(\ul{z}),N_k(\ul{z})}^{Z_k(\ul{z})}\not=0$ for
some finite field $k=\bbF_q$ and $\ul{z}\in\cX_k(\ul{d})$, then both $\bfaz$ and $\bfbz$ are of torsion type.
Thus, the existence of $\vphi_{\bfaz,\bfbz}^{\bfgz}$ follows from Proposition \ref{exist-regular}.

Now assume $Z_k(\ul{z})=Z_k(\ul{z})_t\oplus Z_k(\ul{z})_f$ with $Z_k(\ul{z})_t\neq 0$ and $Z_k(\ul{z})_f\neq 0$
for some finite field $k=\bbF_q$ and $\ul{z}\in\cX_k(\ul{d})$.
Since $\Ext^1(Z_k(\ul{z})_f, Z_k(\ul{z})_t)=0$, applying the formula
\eqref{simple-Green-formula} to the quadruple $(M,N,Z_f,Z_t)=(M_k(\ul{z}),N_k(\ul{z}),Z_k(\ul{z})_f, Z_k(\ul{z})_t)$
gives the equality
$$F_{M,N}^Z = \sum_{A,B,C,D}q^{\langle Z_f,Z_t\rangle-\langle A,D\rangle}
F_{A,B}^MF_{C,D}^NF_{A,C}^{Z_f}F_{B,D}^{Z_t}
\frac{a_Aa_Ba_Ca_D}{a_Ma_N}.$$
 On the one hand, $F_{B,D}^{Z_t}\neq 0$ implies that $B,D\in\cohz\bbX_k$, i.e., $B=B_t$ and $D=D_t$,
 while $F_{A,C}^{Z_f}\neq 0$ implies
$C\in\vect\bbX_k$, i.e., $C=C_f$. On the other hand, since $\Ext^1(\vect\bbX_k, \cohz\bbX_k)=0$,
we obtain that $D=N_t$ and $C=N_f$. Moreover, the associativity of Hall numbers implies that
$$F_{A,B_t}^M = \sum_EF_{A_f,A_t}^EF_{E,B_t}^M = \sum_EF_{A_f,E}^MF_{A_t,B_t}^E = \sum_{E_t}F_{A_f,E_t}^MF_{A_t,B_t}^{E_t} = \delta_{A_f,M_f}F_{A_t,B_t}^{M_t}.$$
 Thus,
$$F_{M,N}^Z=\sum_{A,B\in\cohz\bbX_k}q^{\langle Z_f,Z_t\rangle-\langle A\oplus M_f,N_t\rangle}
F_{A,B}^{M_t}F_{A\oplus M_f, N_f}^{Z_f}F_{B,N_t}^{Z_t}\frac{a_{A\oplus M_f}a_{B}a_{N_f}a_{N_t}}{a_Ma_N}.$$
 Since $M_t=S_k(\bfaz_t,\ul{z})$ with $\bfaz_t=(\az_t,\lz)$, $F_{A,B}^{M_t}\not=0$ implies that
$$A=S_k(\boldsymbol{\xi},\ul{z})\;\text{ and }\; B=S_k(\boldsymbol{\eta},\ul{z})$$
 for $\boldsymbol{\xi}, \boldsymbol{\eta}\in\cS$ of type $\ul{d}$. Moreover, by Proposition \ref{exist-regular},
 there are only finitely many such pairs $(\boldsymbol{\xi}, \boldsymbol{\eta})$. Hence, we conclude that
$$\aligned F_{M,N}^Z=&\sum_{\boldsymbol{\xi}, \boldsymbol{\eta}\in\cS}
q^{\langle Z_f,Z_t\rangle-\langle {S_k(\boldsymbol{\xi},\ul{z})}\oplus M_f,N_t\rangle}
F_{{S_k(\boldsymbol{\xi},\ul{z})},{S_k(\boldsymbol{\eta},\ul{z})}}^{M_t}
F_{{S_k(\boldsymbol{\xi},\ul{z})}\oplus M_f, N_f}^{Z_f}F_{{S_k(\boldsymbol{\eta},\ul{z})},N_t}^{Z_t}\\
&\qquad\quad\times \frac{a_{S_k(\boldsymbol{\xi},\ul{z})\oplus M_f}a_{S_k(\boldsymbol{\eta},\ul{z})}a_{N_f}a_{N_t}}{a_Ma_N},
\endaligned$$
 where the sum, as indicated above, is essentially a finite sum. Applying Proposition \ref{exist-regular}
again shows that $F_{{S_k(\boldsymbol{\xi},\ul{z})},{S_k(\boldsymbol{\eta},\ul{z})}}^{M_t}$
and $F_{{S_k(\boldsymbol{\eta},\ul{z})},N_t}^{Z_t}$ are given by Hall polynomials $\vphi^{\bfaz_t}_{\boldsymbol{\xi},\boldsymbol{\eta}}$
and $\vphi^{\bfgz_t}_{\boldsymbol{\eta},\bfbz_t}$, respectively.  By Lemmas \ref{Ringel dividing lemma} and
\ref{poly for auto}, the existence of Hall polynomial $\vphi_{\bfaz,\bfbz}^{\bfgz}$ follows from that of
Hall polynomials $\vphi^{\bfgz_f}_{\boldsymbol{\xi}\oplus\bfaz_f,\bfbz_f}$.
\end{pf}

\subsection{Reduction 2} To prove Theorem
\ref{existence-HallPoly}, it suffices to prove the existence of Hall
polynomials $\vphi_{\bfaz,\bfbz}^{\bfgz}$ for all $\bfaz\in\cS_t$
and $\bfbz, \bfgz\in \cS_f$.

\begin{pf} By Reduction 1, we can assume that $\bfgz\in \cS_f$ and, thus, $\bfbz\in \cS_f$
since $\vect\bbX$ is closed under subobjects. Thus,
$\bfgz=(\gz,\emptyset)$ and $\bfbz=(\bz,\emptyset)$ for some
$\gz,\bz\in\chi_f$. If $\bfaz\in\cS_f$, then the existence of
$\vphi_{\bfaz,\bfbz}^{\bfgz}$ follows from Proposition
\ref{extistence-vb}.

Now assume $M=S_k(\bfaz,\ul{z})=M_t\oplus M_f$ with $M_t\neq 0$ and $M_f\neq 0$ for some finite
field $k$ and $\ul{z}\in\cX_k(\ul{d})$. Associativity of
Hall numbers implies that
$$F_{M,N}^Z= \sum_EF_{M_f,M_t}^EF_{E,N}^Z = \sum_E F_{M_f,E}^Z F_{M_t, N}^E,$$
where $N=S_k(\bz)$ and $Z=S_k(\gz)$. The term $F_{M_f,E}^Z F_{M_t,
N}^E\not=0$ implies that $E\in\vect\bbX_k$ and there are embeddings
$N\hookrightarrow E \hookrightarrow Z$. Thus, there are only
finitely many isoclasses of such $E$'s and $E=S_k(\thz)$ for some
$\thz\in\cS_f$. Hence,
$$F_{M,N}^Z= \sum_{\thz\in\chi_f} F_{M_f,S_k(\thz)}^Z F_{M_t, N}^{S_k(\thz)}.$$
 By Proposition \ref{extistence-vb}, $F_{M_f,S_k(\thz)}^Z$ is given by the Hall polynomial
$\vphi_{\bfaz_f,\thz}^{\bfgz}$. Therefore, the existence of $\vphi^{\bfgz}_{\bfaz,\bfbz}$ follows
from that of the $\vphi_{\bfaz_t,\bfbz}^\thz$.
\end{pf}

\subsection{Reduction 3} Let $t$ be an integer with $t\geq 2$.
Suppose that Hall polynomials exist for all $\bfaz\in\cS_t$ and
$\bfbz,\bfgz\in\cS_f$ with $\rk \bfgz<t$. Then the Hall polynomial
$\vphi_{\bfaz,\bfbz}^{\bfgz}$ also exists for $\bfaz\in\cS_t$ and
$\bfbz,\bfgz\in\cS_f$ with $\rk \bfgz=t$.

\begin{pf} Take $\bfaz\in\cS_t$ and $\bfbz,\bfgz\in\cS_f$ with $\rk \bfgz=t$.
Let $k=\bbF_q$ be a finite field. If $Z=S_k(\bfgz)$ is decomposable, we can write $\bfgz=\bfgz_1\oplus \bfgz_2$
with $Z_1=S_k(\bfgz_1)\not=0$, $Z_2=S_k(\bfgz_2)\not=0$, and $\Ext^1(Z_2, Z_1)=0$. Applying \eqref{simple-Green-formula}
to the quadruple $(M=S_k(\bfaz,\ul{z}), N=S_k(\bfbz), Z_2,Z_1)$, we obtain the equality
$$F_{M,N}^Z=\sum_{A,B,C,D}q^{\langle Z_2,Z_1\rangle-\langle A,D\rangle}
F_{A,B}^MF_{C,D}^NF_{A,C}^{Z_2}F_{B,D}^{Z_1}\frac{a_Aa_Ba_Ca_D}{a_Ma_N}.$$
 On the one hand, $F_{A,C}^{Z_2}\neq 0$ implies that $C$ is a subobject of $Z_2$ and
 $$\deg C=\deg Z_2-\deg A\geq \deg Z_2-\deg M.$$
Thus, there are only finitely many such $C$'s and each of them has
the form $C=S_k(\thz)$ for some $\thz\in\cS_f$. Similarly,
$F_{B,D}^{Z_1}\neq 0$ implies that $D=S_k(\sz)$ for finitely many
choices of $\sz\in\cS_f$. On the other hand, since
$M=S_k(\bfaz,\ul{z})$ is a torsion sheaf, we have by Proposition
\ref{exist-regular} that $F_{A,B}^M\neq 0$ implies that
$A=S_k(\bfxi,\ul{z})$ and $B=S_k(\bfeta,\ul{z})$ for some
$\bfxi,\bfeta\in\cS_t$ of type $\ul{d}$ and, moreover, there are
only finitely many such pairs $(\bfxi,\bfeta)$. Therefore,
$$\aligned F_{M,N}^Z=\sum_{\bfxi,\bfeta\in\cS_t;\,\thz,\sz\in\cS_f}& q^{\langle Z_2,Z_1\rangle-\langle S_k(\bfxi,\ul{z}),S_k(\sz)\rangle} F_{S_k(\bfxi,\ul{z}),S_k(\bfeta,\ul{z})}^M F_{S_k(\thz),S_k(\sz)}^N F_{S_k(\bfxi,\ul{z}),S_k(\thz)}^{Z_2}\\
& \qquad \times F_{S_k(\bfeta,\ul{z}),S_k(\sz)}^{Z_1}
\frac{a_{S_k(\bfxi,\ul{z})}a_{S_k(\bfxi,\ul{z})}a_{S_k(\thz)}a_{S_k(\sz)}}{a_Ma_N}.
\endaligned$$
 By Propositions \ref{extistence-vb} and \ref{exist-regular}, $F_{S_k(\dz),S_k(\sz)}^N$ and $F_{S_k(\bfxi,\ul{z}),S_k(\bfeta,\ul{z})}^M$
are given by Hall polynomials $\vphi_{\dz,\sz}^{\bfgz}$ and $\vphi^{\bfaz}_{\bfxi,\bfeta}$, respectively.
Since $\rk Z_1<\rk Z=t$ and $\rk Z_2<\rk Z=t$, we have by the induction hypothesis that $F_{S_k(\bfxi,\ul{z}),S_k(\thz)}^{Z_2}$
and $F_{S_k(\bfeta,\ul{z}),S_k(\sz)}^{Z_1}$ are given by Hall polynomials, too. Applying
Lemmas \ref{poly for auto} and \ref{Ringel dividing lemma} gives the existence of $\vphi_{\bfaz,\bfbz}^{\bfgz}$.

Now suppose that $Z=S_k(\bfgz)$ is indecomposable.
Since $\rk Z=t\geq 2$, it follows from Lemma \ref{vector bundle factorization} that there is an exact sequence
$0\to L\to Z\to Z'\to 0$ with $L$ a line bundle and $Z'$ a vector
bundle such that $\Ext^1(Z', L)\cong k$. Applying \eqref{Green-formula} to the quadruple $(M,N,Z',L)$,
we have
$$\sum_EF_{M,N}^EF_{Z', L}^E/a_E =\sum_{A,B,C,D}q^{-\langle A,D\rangle}F_{A,B}^MF_{C,D}^NF_{A,C}^{Z'}F_{B,D}^{L}\frac{a_Aa_Ba_Ca_D}{a_Ma_Na_{Z'}a_{L}}.$$
 Since every extension of $Z'$ by $L$ is isomorphic either to $Z$ or to $Z'\oplus L$, it follows that
$$\sum_EF_{M,N}^EF_{Z', L}^E/a_E = F_{M,N}^Z F_{Z',L}^Z/a_Z+F_{M,N}^{Z'\oplus L}F_{Z',L}^{Z'\oplus L}/a_{Z'\oplus L}.$$
 Consequently, we obtain that
$$\aligned
F_{M,N}^Z=&\frac{a_Z}{F_{Z',L}^Z}\sum_{A,B,C,D}q^{-\langle A,D\rangle}F_{A,B}^MF_{C,D}^NF_{A,C}^{Z'}F_{B,D}^{L}
\frac{a_Aa_Ba_Ca_D}{a_Ma_Na_{Z'}a_{L}}\\
&-\frac{a_Z}{a_{Z'\oplus L} F_{Z',L}^Z}F_{M,N}^{Z'\oplus L}F_{Z',L}^{Z'\oplus L}.
\endaligned$$
By similar arguments as above, the sum on the right-hand side is
finite and each Hall number occurring in the sum, as well as $F_{M,N}^{Z'\oplus L}$
(since $Z'\oplus L$ is decomposable), is given by an integer polynomial. By Lemma
\ref{Riedtmann-Peng-formula}, it is direct to see that
$F_{Z',L}^{Z'\oplus L}$ is given by an integer polynomial. Further,
by Proposition \ref{extistence-vb}, $F_{Z',L}^Z$ is given by a monic
polynomial. Therefore, $\vphi_{\bfaz,\bfbz}^{\bfgz}$ exists by Lemma
\ref{Ringel dividing lemma}.
\end{pf}

\subsection{The proof of Theorem \ref{existence-HallPoly}}
Combining Reductions 1, 2 and 3, we are reduced to prove the
existence of the Hall polynomials $\vphi_{\bfaz,\bfbz}^{\bfgz}$ for
the case where $\bfaz\in \cS_t$, $\bfbz, \bfgz\in \cS_f$ with
$\rk\bfgz=1$.

Let $k$ be a field and $\ul{z}=(z_1,\ldots,z_r)\in\cX_k(\ul{d})$. As above, put
$Z=S_k(\bfgz)$, $M=S_k(\bfaz,\ul{z})$ and $N=S_k(\bfbz)$. By taking a grading shift, we may assume
that $N=\co$ and $Z=\co(\vu)$ for some $\vu\in\mathbb{L}_+$. We now use associativity together
with an induction on the determinant $\det\vu$ to prove the assertion.

If $M$ supports at more than two distinct points, then $M$ admits a
non-trivial decomposition $M=M_1\oplus M_2$ such that $M_1$ and
$M_2$ have disjoint supports. It follows that
$$\Ext^1(M_1,M_2)=0=\Ext^1(M_2,M_1).$$ By associativity, we have
$$F_{M,\co}^{\co(\vu)}= \sum_EF_{M_2,M_1}^EF_{E, \co}^{\co(\vu)}=\sum_{\co(\vv)}F_{M_2,\co(\vv)}^{\co(\vu)}F_{M_1,\co}^{\co(\vv)}
=\sum_{0<\vv<\vu}F_{M_2(-\vv),\co}^{\co(\vu-\vv)}F_{M_1,\co}^{\co(\vv)}.$$
By induction on $\det \vu$, all the terms on the right-hand side are
given by Hall polynomials. This shows the existence of
$\vphi_{\bfaz,\bfbz}^{\bfgz}$.

Now we assume that $M$ supports at a single point. Then $F_{M,\co}^{\co(\vu)}\neq 0$
implies that there is a surjection $\co(\vu)\twoheadrightarrow M$, which ensures that $M$ is
indecomposable; see, e.g., \cite[Ex.~4.12]{Sch2006}. If $M$ has quasi-Loewy length $\ell\geq 2$, then
there is an exact sequence $0\to M'\to M\to S\to 0$, where $S$ is
quasi-simple. Associativity of Hall numbers implies that
$$F_{S, M'}^{M}F_{M, \co}^{\co(\vu)}+F_{S, M'}^{S\oplus M'}F_{S\oplus M',\co}^{\co(\vu)}
=\sum_{\co(\vv)}F_{S,\co(\vv)}^{\co(\vu)}F_{M',\co}^{\co(\vv)}
=\sum_{0<\vv<\vu}F_{S(-\vv),\co}^{\co(\vu-\vv)}F_{M',\co}^{\co(\vv)}.$$
 Clearly, $F_{S\oplus M',\co}^{\co(\vu)}=0$ and $F_{S, M'}^{M}=1$. Thus,
$$F_{M, \co}^{\co(\vu)}=\sum_{0<\vv<\vu}F_{S(-\vv),\co}^{\co(\vu-\vv)}F_{M',\co}^{\co(\vv)}$$
 is given by an integer polynomial by an inductive argument on $\det \vu$.

 Now assume $M$ is a quasi-simple sheaf. If $M=S_{i,j}$ lies in a non-homogeneous tube,
then $F_{S_{i,j},\co}^{\co(\vu)}=1$ for $\vu=\vx_i$ and $j=1$, and
zero otherwise. If $M\in{\rm coh}_{z}\text{-}\bbX_k$ for some $z\in\bbH_k$ with $r=\deg(z)$,
then $F_{M, \co}^{\co(\vu)}=1$ for $\vu=r\vc$ and it is zero
otherwise. In both cases, the Hall polynomial
$\vphi_{\bfaz,\bfbz}^{\bfgz}$ exists.

This finishes the proof of Theorem \ref{existence-HallPoly}.

\medskip

Let $k$ be a finite field and $X\in\coh\bbX_k$. Then for each field extension $k\subseteq K$,
$X^K:=X\otimes_kK$ is an object in $\coh\bbX_K$. A finite field extension
$K$ of $k$ is said to be {\it conservative} relative to $X$ if for each indecomposable summand
$Y$ of $X$, $Y^K$ is indecomposable in $\coh\bbX_K$. In general, given a finite
set ${\scr X}=\{X_1,\ldots,X_m\}$ of objects in $\coh\bbX_k$, a finite field extension
$K$ of $k$ is said to be {\it conservative} relative to $\scr X$ if $K$ is conservative relative
to each $X_i$ for $1\leq i\leq m$. Note that there always exist infinitely many conservative field extensions
of $k$ relative to $\scr X$.

\begin{Cor} \label{Hall-poly-3-obj} Fix a finite field $k$ and three objects $M,N,Z$ in $\coh\bbX_k$.
Then there exists a polynomial $\vphi_{M,N}^Z\in\bbZ[T]$ such that for each conservative field extension $K$
of $k$ relative to $\{M,N,Z\}$,
$$\vphi_{M,N}^Z(|K|)=F_{M^K,N^K}^{Z^K}.$$
\end{Cor}

\begin{pf} Choose decomposition sequences $\bfaz=(\az,\lz),\bfbz=(\bz,\mu),\bfgz=(\gz,\nu)$
of the same type, say of type $\ul{d}=(d_1,\ldots,d_r)$, and $\ul{z}=(z_1,\ldots,z_r)\in \cX_k(\ul{d})$ such that
$$M=S_k(\bfaz,\ul{z}),\;\,N=S_k(\bfbz,\ul{z}),\;\text{ and }\;Z=S_k(\bfgz,\ul{z}).$$
 Moreover, we can assume that for each $1\leq i\leq r$, one of the partitions $\lz^{(i)},\mu^{(i)},\nu^{(i)}$
 is not the empty partition. Thus, if $K$ is a conservative field extension
of $k$ relative to $\{M,N,Z\}$, then $\ul{z}\in\cX_K(\ul{d})$ and
$$M^K=S_K(\bfaz,\ul{z}),\;\,N^K=S_K(\bfbz,\ul{z}),\;\text{ and }\;Z^K=S_K(\bfgz,\ul{z}).$$
 By Theorem \ref{existence-HallPoly}, $\vphi^{\bfgz}_{\bfaz,\bfbz}$ is the desired
 polynomial $\vphi_{M,N}^Z$.
\end{pf}

\section{Generic Hall algebra of $\bbX$ and its Drinfeld double}

In this section, we define the (generic) Hall algebra of a domestic weighted projective line
$\bbX$ as well as its Drinfeld double by
using the Hall polynomials given in the previous sections.

\subsection{Generic Hall algebra} By 2.2, for each finite field $k$,
we have the Ringel--Hall algebra $H(\bbX_k)$ of the category of coherent
sheaves on $\bbX_k$ defined over $k$.

Recall the set of decomposition classes $\cS$ over $\bbX$ defined in Section 3 and
the Hall polynomial $\vphi^{\bfgz}_{\bfaz,\bfbz}(T)\in \bbZ[T]$ for each triple $\bfaz, \bfbz,\bfgz\in\cS$ of the same type.
Let $\bbQ[\bfv,\bfv^{-1}]$ be the Laurent polynomial ring with indeterminate $\bfv$ and
put
$$H_\bfv(\bbX):=\bigoplus\limits_{\bfaz\in\cS}\bbQ[\bfv, \bfv^{-1}]u_{\bfaz},$$
that is, the free $\bbQ[\bfv,\bfv^{-1}]$-module with basis $\{u_{\bfaz}\mid\bfaz\in\cS\}$.
For $\bfaz, \bfbz\in\cS$, define their multiplication by
$$u_{\bfaz}u_{\bfbz}=\bfv^{\langle\bfaz, \bfbz\rangle}\sum_{\bfgz\in\cS}\vphi_{\bfaz,\bfbz}^{\bfgz}(\bfv^2)u_{\bfgz},$$
 where $\bfaz$ and $\bfbz$ are thought of same type in the sense of
Remark \ref{equiv-Segre-class} and the sum is taken over all $\bfgz\in\cS$ of the same type.
Note that for fixed $\bfaz,\bfbz\in\cS$, there are only finitely many $\bfgz$ satisfying
$\vphi_{\bfaz,\bfbz}^{\bfgz}(T)\not=0$. If $\bfaz=(\az,\emptyset)$, we sometimes
write $u_{\bfaz}=u_{[S(\az)]}$ for computational purpose, e.g., $u_{[\cal O]}$, $u_{[S_{i,j}]}$, etc.

\begin{Prop} The $\bbQ[\bfv,\bfv^{-1}]$-module $H_\bfv(\bbX)$ endowed with the multiplication defined above
becomes an associative algebra with identity $1=u_0$, where $0$ denotes
the decomposition class $(0,\emptyset)$.
\end{Prop}

\begin{pf} We need to show the associativity of the multiplication.
Take arbitrary $\bfaz,\bfbz,\bfgz\in\cS$. On the one hand, we have
$$\aligned
(u_{\bfaz}u_{\bfbz})u_{\bfgz}&=\bfv^{\langle\bfaz,\bfbz\rangle}\sum_{\bfthz\in\cS}\vphi_{\bfaz,\bfbz}^{\bfthz}(\bfv^2)u_{\bfthz}u_{\bfgz}
=\bfv^{\langle\bfaz,\bfbz\rangle}\sum_{\bfthz\in\cS}\vphi_{\bfaz,\bfbz}^{\bfthz}(\bfv^2)
\big(\bfv^{\langle\bfthz,\bfgz\rangle}\sum_{\bfdz\in\cS}\vphi_{\bfthz,\bfgz}^{\bfdz}(\bfv^2)u_{\bfdz}\big)\\
&=\bfv^{\langle\bfaz,\bfbz\rangle+\langle\bfaz,\bfgz\rangle+\langle\bfbz,\bfgz\rangle}
\sum\limits_{\bfdz\in\cS}\big(\sum_{\bfthz\in\cS}\vphi_{\bfaz,\bfbz}^{\bfthz}(\bfv^2)\vphi_{\bfthz,\bfgz}^{\bfdz}(\bfv^2)\big)u_{\bfdz}.
\endaligned$$
On the other hand,
$$u_{\bfaz}(u_{\bfbz}u_{\bfgz})=\bfv^{\langle\bfaz,\bfbz\rangle+\langle\bfaz,\bfgz\rangle+\langle\bfbz,\bfgz\rangle}
\sum\limits_{\bfdz\in\cS}\big(\sum_{\bfthz\in\cS}\vphi_{\bfaz,\bfthz}^{\bfdz}(\bfv^2)\vphi_{\bfbz,\bfgz}^{\bfthz}(\bfv^2)\big)u_{\bfdz}.$$
Thus, to prove the associativity, it suffices to show that
\begin{equation}\label{asso}
\sum_{\bfthz\in\cS}\vphi_{\bfaz,\bfbz}^{\bfthz}(T)\vphi_{\bfthz,\bfgz}^{\bfdz}(T)
=\sum_{\bfthz\in\cS}\vphi_{\bfaz,\bfthz}^{\bfdz}(T)\vphi_{\bfbz,\bfgz}^{\bfthz}(T).
\end{equation}
By the definition, for each finite field $k$ with $q=|k|\gg0$ and $\ul{z} \in\cX_{k}(\ul{d})$,
$$\aligned
\sum_{\bfthz\in\cS}\vphi_{\bfaz,\bfbz}^{\bfthz}(q)\vphi_{\bfthz,\bfgz}^{\bfdz}(q)
&=\sum_{\bfthz\in\cS}F_{S_k(\bfaz,\ul{z}), S_k(\bfbz,\ul{z})}^{S_k(\bfthz,\ul{z})} F_{S_k(\bfthz,\ul{z}),S_k(\bfgz,\ul{z})}^{S_k(\bfdz,\ul{z})}\\ &=F_{S_k(\bfaz,\ul{z}),S_k(\bfbz,\ul{z}),S_k(\bfgz,\ul{z})}^{S_k(\bfdz,\ul{z})}
=\sum_{\bfthz\in\cS}\phi_{\bfaz,\bfthz}^{\bfdz}(q)\phi_{\bfbz,\bfgz}^{\bfthz}(q),
\endaligned$$
where the second equality follows from the associativity of the Ringel--Hall algebra
$H(\bbX_k)$. Hence, the left and right hand sides of (\ref{asso}) take the same values
for prime power $q\gg0$. In conclusion, (\ref{asso}) holds.
\end{pf}

Form now onwards, let $\frakQ$ denote the set of all prime powers ($\neq 1$). For each $q\in\frakQ$,
let $\bbF_q$ denote the field with $q$ elements and set $v_q=\sqrt{q}$. We will simply write $\bbX_q=\bbX_{\bbF_q}$,
$S_q(\bfaz,\ul{z})=S_{\bbF_q}(\bfaz,\ul{z})$, etc.  Consider the infinite direct product
$$\prod_{q\in\frakQ}H(\bbX_q)$$
which clearly carries a natural structure of an associative algebra whose multiplication is defined
componentwise, that is, for any $(a_q)_q, (b_q)_q\in\prod_q H(\bbX_q)$,
$$(a_q)_q\cdot(b_q)_q=(a_q b_q)_q.$$
Let $\mathcal {I}$ be the ideal of $\prod_{q}H(\bbX_q)$ generated by the elements $(a_q)_q$
with $a_q=0$ for $q\gg0$. Then the quotient
$$\widehat{H(\bbX)}=:\prod_{q\in\frakQ}H(\bbX_q)/\mathcal {I}$$
becomes an associative algebra, too. In other words, elements in $\widehat{H(\bbX)}$
are equivalence classes under the equivalence relation $\sim$ defined by
$$(a_q)_q\sim (b_q)_q\Longleftrightarrow a_q=b_q\;\text{ for $q\gg0$.}$$
 It is clear that the element $\wt v=(v_q)_q$ is invertible and does not satisfy any polynomial
equation over $\bbQ$ in $\widehat{H(\bbX)}$. Thus, by identifying $\wt v$ with $\bfv$, $\widehat{H(\bbX)}$
can be viewed as an algebra over the Laurent polynomial ring $\bbQ[\bfv,\bfv^{-1}]$.

In what follows, for each $q\in\frakQ$, we fix a total ordering $\preccurlyeq$ of all points in $\bbH_{\bbF_q}$
such that $x\preccurlyeq y$ implies $\deg(x)\leq\deg(y)$. The chain of points in $\bbH_{\bbF_q}$ of degree $d$
is denoted by
$$y_{d,1}\prec y_{d,2}\prec\cdots\prec y_{d,\zeta_d},$$
 where $\zeta_d$ is the number of points of degree $d$ in $\bbH_{\bbF_q}$. For each given $\bfaz\in\cS$ of
type $\ul{d}=(d_1,\ldots,d_r)$, we fix an element $\ul{z}=\ul{z}_{\bfaz,q}=(z_1,\ldots,z_r)\in\cX_{\bbF_q}(\ul{d})$
such that for each $d_1\leq d\leq d_r$,
$$(z_i,\ldots,z_j)=(y_{d,1},\ldots, y_{d,j-i+1}),$$
 where $d_{i-1}<d=d_{i}=\cdots=d_j<d_{j+1}$. In fact, $\ul{z}$ is independent of $\bfaz$.
In particular, if two decomposition sequences $\bfaz=(\az,\lz)$ and $\bfbz=(\bz,\mu)$ can be identified in
the sense of Remark \ref{equiv-Segre-class},
then $S_q(\bfaz,\ul{z}_{\bfaz,q})\cong S_q(\bfbz,\ul{z}_{\bfbz,q})$.

\begin{Prop} \label{embedding-Phi} The assignment
$u_{\bfaz}\mapsto ([S_q(\bfaz,\ul{z}_{\bfaz,q})])_q, \bfaz\in\cS$, defines an embedding of
$\bbQ[\bfv,\bfv^{-1}]$-algebras
$$\Phi:H_\bfv(\bbX)\lra\widehat{H(\bbX)}=\prod_{q\in\frakQ} H(\bbX_q)/\mathcal {I}.$$
\end{Prop}

\begin{pf} The injectivity of $\Phi$ is obvious. We need to check that $\Phi$ is
an algebra homomorphism. For any $\bfaz,\bfbz\in\cS$ (of same type), we have by the definition that
$$u_{\bfaz}u_{\bfbz}=\bfv^{\langle\bfaz,\bfbz\rangle}\sum_{\bfgz\in\cS}\vphi_{\bfaz,\bfbz}^{\bfgz}(\bfv^2)u_{\bfgz}.$$
Furthermore, the $q$-component of $\Phi(u_{\bfaz})\Phi(u_{\bfbz})$ is given by
$$[S_q(\bfaz,\ul{z})][S_q(\bfbz,\ul{z})]=
v_q^{\langle S_q(\bfaz,\ul{z}),S_q(\bfbz,\ul{z})\rangle}\sum\limits_{S_q(\bfgz,\ul{z})}
F_{S_q(\bfaz,\ul{z}),S_q(\bfbz,\ul{z})}^{S_q(\bfgz,\ul{z})}[S_q(\bfgz,\ul{z})],$$
 where $\ul{z}=\ul{z}_{\bfaz,q}=\ul{z}_{\bfbz,q}=\ul{z}_{\bfgz,q}$,
while the $q$-component of $\Phi(u_{\bfaz}u_{\bfbz})$ is given by
$$v_q^{\langle\bfaz,\bfbz\rangle}\sum_{\bfgz}\vphi_{\bfaz,\bfbz}^{\bfgz}(q)[S_q(\bfgz,\ul{z})].$$
By Theorem \ref{existence-HallPoly},
$\vphi_{\bfaz,\bfbz}^{\bfgz}(q)=F_{S_q(\bfaz,\ul{z}),S_q(\bfbz,\ul{z})}^{S_q(\bfgz,\ul{z})}$ for $q\gg0$.
Hence, $\Phi(u_{\bfaz}u_{\bfbz})=\Phi(u_{\bfaz})\Phi(u_{\bfbz})$, as desired.
\end{pf}

From now on, we will identify $H_\bfv(\bbX)$ with the subalgebra $\Phi(H_\bfv(\bbX))$
of $\widehat{H(\bbX)}$. Thus, we will use the notation $u_{\bfaz}$ to denote its image
$\Phi(u_{\bfaz})$ in $\widehat{H(\bbX)}$.

Recall the elements $T_{r,q}$, $Z_{r,q}$ and $\Thz_{\vx,q}$ in $H(\bbX_q)$ introduced in \S 2.3
for all $r\geq 1$, $\vx\in\bbL$, and $q\in\frakQ$. Set
$$Z_r:=(Z_{r,q})_q,\; T_r:=(T_{r,q})_q,\;\Thz_{\vx}:=(\Thz_{\vx,q})_q\in \widehat{H(\bbX)}=\prod_{q\in\frakQ}H(\bbX_q)/\cI.$$
 Let $\calH_\bfv(\bbX)$ be the $\bbQ[\bfv,\bfv^{-1}]$-subalgebra
 of $\widehat{H(\bbX)}$ generated by $H_\bfv(\bbX)$ and $Z_r$ for $r\geq 1$ and put
 $$\bfcalH_\bfv(\bbX)=\calH_\bfv(\bbX)\otimes_{\bbQ[\bfv,\bfv^{-1}]}\bbQ(\bfv).$$
Both $\calH_\bfv(\bbX)$ and $\bfcalH_\bfv(\bbX)$ are called the {\it generic Hall algebras} of $\bbX$.
%By the definition, $\bfcalH_\bfv(\bbX)$ is spanned by the set
%\begin{equation} \label{spann-set-generic-HA}
%\{ u_{\bfaz}\prod_{r\geq 1}Z_r^{l_r}\mid \bfaz\in\cS, l_r\geq 0\;\text{with}\; \sum_{r\geq 1}l_r<\infty \}.
%\end{equation}

\begin{Prop}\label{common subalgebra} For all $r\geq 1$ and $\vx\in\bbL_+$, $T_r$ and $\Thz_{\vx}$ lie in
$\bfcalH_\bfv(\bbX)$. Moreover, $\bfcalH_\bfv(\bbX)$ can be also generated by $H_\bfv(\bbX)$
together with one of the following sets:
$${\rm (I)}\;\, \{T_r\mid  r\geq 1\};\;\; {\rm (II)}\;\, \{\Thz_{r\vc}\mid r\geq 1\};\;\;
{\rm (III)}\;\, \{\Thz_{\vx}\mid \vx\in\bbL_+\}.$$
\end{Prop}

\begin{pf} By \cite[Prop.~5.6]{BurSch}, for every $q\in\frakQ$ and $r>0$, $T_{r,q}-Z_{r,q}$ belongs
to the subalgebra of $H(\bbX_q)$ generated by $[S_{i,j}]$ for $i\in I$ and $0\leq
j\leq p_i-1$. It follows
that all $T_r-Z_{r}$ belong to the subalgebra of $\bfcalH_\bfv(\bbX)$ generated by the $u_{[S_{i,j}]}$.
Hence,
$$\langle H_\bfv(\bbX), T_r\mid r\geq 1\rangle=\langle H_\bfv(\bbX), Z_r\mid r\geq 1\rangle.$$

By \cite[Ex.~4.12]{Sch2006} and \cite[Lem.~5.20]{BurSch}, we get that
\begin{equation}\label{relation of Tr and Thz_r}
1+\sum_{r\geq 1}\Thz_{r\vc}t^r={\rm exp}\big((\bfv-\bfv^{-1})\sum_{r\geq 1}T_{r}t^r\big).
\end{equation}
This implies that
$$\langle H_\bfv(\bbX), T_r\mid r\geq 1\rangle=\langle H_\bfv(\bbX), \Thz_{r\vc}\mid r\geq 1\rangle.$$
Moreover, by \cite[Prop.~5.21]{BurSch}, we have
$\Thz_{\vx}\in \langle H_\bfv(\bbX), \Thz_{r\vc}\mid r\geq 1\rangle$ for $\vx\in\bbL_+$.
This finishes the proof.
\end{pf}

Let $\mathcal{L}$ denote the set of infinite sequences of nonnegative integers $\ul{l}=(l_r)_{r\geq 1}$
satisfying $\sum_r l_r<\infty$. For each $\ul{l}\in\mathcal{L}$, set
\[Z_{\ul{l}}=\prod\limits_{r\geq 1}Z_r^{l_r}, \quad T_{\ul{l}}=\prod\limits_{r\geq 1}T_r^{l_r}, \text{\quad and\quad} \Thz_{\ul{l}}=\prod\limits_{r\geq 1}\Thz_{r\vc}^{l_r}.
\]

\begin{Prop}\label{three-bases-generic-HA} Each of the following three sets
$$\{u_{\bfaz}Z_{\ul{l}}\mid \bfaz\in\cS, \ul{l}\in{\cal L}\},
\quad \{u_{\bfaz}T_{\ul{l}}\mid \bfaz\in\cS, \ul{l}\in{\cal L}\},
\quad \text{and} \quad \{u_{\bfaz}\Thz_{\ul{l}}\mid \bfaz\in\cS, \ul{l}\in\mathcal{L}\} $$
is a $\bbQ(\bfv)$-basis of $\bfcalH_\bfv(\bbX)$.
\end{Prop}

\begin{pf} We only prove that ${\scr Z}:=\{u_{\bfaz}Z_{\ul{l}}\mid \bfaz\in\cS, \ul{l}\in{\cal L}\}$
is a $\bbQ(\bfv)$-basis of $\bfcalH_\bfv(\bbX)$. This together with the arguments in the proof of Proposition
\ref{common subalgebra} implies that the other two sets are $\bbQ(\bfv)$-bases, too.

By the definition, $\bfcalH_\bfv(\bbX)$ is spanned by ${\scr Z}$. It remains to show that $\scr Z$
is a linearly independent set. For each $\bfaz\in\cS$ and $\ul{l}\in{\cal L}$, we have
$$u_{\bfaz}Z_{\ul{l}}=u_{\bfaz_f}u_{\bfaz_t}Z_{\ul{l}},$$
where $\bfaz=\bfaz_f\oplus\bfaz_t$ with $\bfaz_f\in\cS_f$ and $\bfaz_t\in\cS_t$.
Then for each $q\in\frakQ$, the $q$-th component of $u_{\bfaz}Z_{\ul{l}}$ is
$$[S_q(\bfaz_f)][S_q(\bfaz_t,\ul{z}_{\bfaz,q})]Z_{\ul{l},q},$$
 where $Z_{\ul{l},q}=\prod_{r\geq 1}(Z_{r,q})^{l_r}$.
Let $H^f(\bbX_q)$ (resp., $H^t(\bbX_q)$) be the $\bbQ[v_q,v_q^{-1}]$-subalgebra of $H(\bbX_q)$ generated
by $[S]$ with $S$ all the torsion-free (resp., torsion) sheaves. Then the multiplication
map
$$H^f(\bbX_q)\otimes_{\bbQ[v_q,v_q^{-1}]}H^t(\bbX_q)\stackrel{{\rm mult}}{\lra} H(\bbX_q)$$
 is an isomorphism of $\bbQ[v_q,v_q^{-1}]$-modules. Thus, it suffices to prove that the set
$$\{[S_q(\bfaz,\ul{z}_{\bfaz,q})]Z_{\ul{l}}\mid \bfaz\in\cS_t, \ul{l}\in{\cal L}\}$$
 is linearly independent in $H(\bbX_q)$ for $q\gg0$. By the construction of $Z_{r,q}$ given in \cite{Sch2006},
 this is reduced to prove that for each fixed $\ul{l}\in{\cal L}$, the set
 $\{[S_q(\bfaz,\ul{z}_{\bfaz,q})]Z_{\ul{l}}\mid \bfaz\in\cS_t\}$ is linearly independent.

 By \cite{Sch2000,Hub2} and \cite[Sect.~5]{BurSch}, the elements $Z_{r,q}$ are
central in $H^t(\bbX_q)$, and for each $\ul{l}\in{\cal L}$, multiplication by
$Z_{\ul{l},q}$
$$H^t(\bbX_q)\lra H^t(\bbX_q),\;a\lmto a Z_{\ul{l},q}$$
 is injective. Therefore, for $q\gg0$, the set
$$\{[S_q(\bfaz,\ul{z})]Z_{\ul{l},q}\mid \bfaz\in\cS_t\}$$
 is linearly independent. We conclude that
 $${\scr Z}=\{u_{\bfaz}Z_{\ul{l}}\mid \bfaz\in\cS, \ul{l}\in{\cal L}\}$$
is a linearly independent set, as desired.
\end{pf}

\begin{Prop}\label{formula of product} The following relations hold in $\bfcalH_\bfv(\bbX)$:
\begin{itemize}
\item[(1)] $[Z_r, Z_s]=0;$
\item[(2)] $[Z_r, u_{\bfaz_{t}}]=0;$
\item[(3)] $[Z_r, u_{[\co(\vx)]}]=\gz_r u_{[\co(\vx+r\vc)]}$ for some $\gz_r\in\bbQ(\bfv)$.
\item[(4)] $[T_r, u_{[\co(\vx)]}]=\frac{[2r]}{r}u_{[\co(\vx+r\vc)]}$, where $[2r]=(\bfv^{2r}-\bfv^{-2r})/(\bfv-\bfv^{-1})$.
\item[(5)] $[\Thz_{r\vc}, u_{[\co(n\vc)]}]=(1-\bfv^{-4})\sum\limits_{1\leq i\leq r}\bfv^{2i}u_{[\co((n+i)\vc)]}\Thz_{(r-i)\vc}.$
\end{itemize}
\end{Prop}

\begin{pf} The relations (1)--(4) hold since they hold in each $q$-component for
$q\in\frakQ$; see \cite{BurSch}.

We now prove (5). For each finite field $k=\bbF_q$, the well-known embedding $\coh \bbP^1(k)\ra \coh\bbX_k$
induces an algebra embedding between their Ringel--Hall algebras which takes
$\Thz_{r,q}\mapsto \Thz_{r\vc,q}$ and $[\co(n)]\mapsto[\co(n\vc)]$; see
\cite[Lem.~5.20]{BurSch}. Hence, it suffices to show that the equality
\begin{equation}  \label{thz-o-at-q}
[\Thz_{r,q}, [\co(n)]]=(1-v^{-4})\sum_{1\leq i\leq r}v^{2i}[\co((n+i))]\Thz_{r-i,q}
\end{equation}
holds in the Ringel--Hall algebra of $\bbP^1(k)$, where $v=\sqrt{q}$. Consider the generating functions
$$\Thz(s)=1+\sum_{i\geq 1}\Thz_{i,q}s^{i},\;\bbT(s)=1+\sum_{i\geq 1}\frac{T_{i}}{[i]_{v}}s^{i},\;\text{ and }\;
\bbO(t)=1+\sum_{i\geq 1}[\co(i)]t^{i},$$
where $[i]_{v}=(v^i-v^{-i})/(v-v^{-1})$. By the proof of \cite[Ex.~4.12]{Sch2006},
$$\Thz(s)=\exp(\bbT(v s)-\bbT(v^{-1}s))\;\text{ and }\;
[\bbT(s), \bbO(t)]=-\bbO(t)\log(1-\frac{s}{v t})(1-\frac{v s}{t}).$$
Hence,
$$[\bbT(v s)-\bbT(v^{-1}s), \bbO(t)]=\bbO(t)\log\frac{1-\frac{s}{v^{2}t}}{1-\frac{v^{2}s}{t}}.$$
This implies that
$$\Thz(s)\bbO(t)=\bbO(t)\Thz(s)\frac{1-\frac{s}{v^{2}t}}{1-\frac{v^{2}s}{t}}.$$
Then \eqref{thz-o-at-q} follows from the equality
$$\frac{1-\frac{s}{v^{2}t}}{1-\frac{v^{2}s}{t}}=1+(1-v^{-4})\sum_{i\geq 1}v^{2i}\frac{s^i}{t^i}.$$
\end{pf}

\subsection{Extended generic Hall algebra}\label{section of ext gen Hall alg}
For each prime power $q\in\frakQ$, let $K_0(\bbX_q)$ be the Grothendieck group
of $\coh\bbX_q$. It is well known that $K_0(\bbX_q)$ is independent
of $q$. Thus, we simply write $K_0(\bbX)$ for $K_0(\bbX_q)$. Finally, set
$$\dz:=[\co(\vec{c})]-[\co]\in K_0(\bbX).$$
 Furthermore, we denote by $\calK=\bbQ[\bfv,\bfv^{-1}][K_0(\bbX)]$
the group algebra of $K_0(\bbX)$. To avoid possible confusion, we write $K_{[M]}$ instead of $[M]$
for each $M\in \coh\bbX$. Furthermore, for each decomposition sequence $\bfaz$ in $\cS$,
we simply set $K_{\bfaz}=K_{[S_q(\bfaz,\ul{z})]}$ in $\calK$ for $q\gg0$. We equip the $\bbQ[\bfv,\bfv^{-1}]$-module
$\wt\calH_\bfv(\bbX):=\calH_\bfv(\bbX)\otimes_{\bbQ[\bfv,\bfv^{-1}]}\calK$
with an algebra structure (containing $\calH_\bfv(\bbX)$ and $\calK$
as subalgebras) by imposing the relations
$$K_{\bfa}u_{\bfaz}K_{\bfa}^{-1}=\bfv^{(\bfa, \bfaz)}u_{\bfaz},\;\forall\,\bfaz\in\cS, \bfa\in K_0(\bbX).$$

For each $q\in\frakQ$, denote by $\wt H(\bbX_q)$ the extended version of
$H(\bbX_q)$, that is,
$$\wt H(\bbX_q)=H(\bbX_q)\otimes_{\bbQ[v_q,v_q^{-1}]}\bbQ[v_q,v_q^{-1}][K_0(\bbX)].$$
Then the embedding of $\bbQ[\bfv,\bfv^{-1}]$-algebras
$$\Phi:H_\bfv(\bbX)\longrightarrow \widehat{H(\bbX)}=\prod_{q\in\frakQ}H(\bbX_q)/\cI$$
extends to an embedding
$$\wt\Phi:\wtH_\bfv(\bbX)\longrightarrow\prod_{q\in\frakQ}\wtH(\bbX_q)/\wt{\cI},\;
u_{\bfaz} K_{\bfa}\longmapsto([S_q({\bfaz},\ul{z})]K_{\bfa})_q,$$
where $\wt\cI$ denotes the ideal of $\prod_{q}\wtH(\bbX_q)$ generated by the elements $(a_q)_q$
with $a_q=0$ for $q\gg0$. We view $\wtH_\bfv(\bbX)$ as a subalgebra of
$\prod_{q}\wtH(\bbX_q)/\wt{\cI}$ and denote by $\wt\calH_\bfv(\bbX)$ the $\bbQ[\bfv,\bfv^{-1}]$-subalgebra
generated by $\wtH_\bfv(\bbX)$ together with $Z_r$ for $r\geq 1$.
Finally, we set
 $$\wt\bfcalH_\bfv(\bbX)=\wt\calH_\bfv(\bbX)\otimes_{\bbQ[\bfv,\bfv^{-1}]}\bbQ(\bfv)
 =\bfcalH_\bfv(\bbX)\otimes_{\bbQ(\bfv)}\bbQ(\bfv)[K_0(\bbX)].$$
  Both $\wt\calH_\bfv(\bbX)$ and $\wt\bfcalH_\bfv(\bbX)$ are called the
\emph{extended generic Hall algebras} of $\bbX$.

The topological comultiplications $\Dz_q$ on the Hall algebra
$\wtH(\bbX_q)$ for all $q\in\frakQ$ induce a topological comultiplication $\Dz$ on
$\prod_{q}\wtH(\bbX_q)$ by $\Dz((a_q)_q)=(\Dz_q(a_q))_q$, where
$a_q\in\wtH(\bbX_q)$. It is easy to see that $\wt{\cI}$ is a coideal. Hence,
it induces a topological comultiplication on the
quotient $\prod_{q}\wtH(\bbX_q)/\wt{\cI}$, which is still denoted
by $\Dz$.

\begin{Prop}\label{comultiplication formula}
We have the following formulas in $\wt\bfcalH_\bfv(\bbX)$:
\begin{itemize}
\item[(1)] $\Dz(u_{\bfgz}K_{\bfa})=\sum\limits_{\bfaz,\bfbz}\bfv^{\langle\bfaz,\bfbz\rangle}\vphi_{\bfaz,\bfbz}^{\bfgz}
\frac{a_{\bfaz} a_{\bfbz}}{a_{\bfgz}}u_{\bfaz}K_{\bfbz+\bfa}\otimes
u_{\bfbz}K_{\bfa}$ for $\bfgz\in\chi_t$ and $\bfa\in K_0(\bbX)$;

\item[(2)] $\Dz(u_{[\co]})=u_{[\co]}\otimes 1+\sum_{\vx\in \bbL_+}\Thz_{\vx}K_{[\co(-\vx)]}\otimes
u_{[\co(-\vx)]}$;

\item[(3)] $\Dz(Z_r)=Z_r\otimes 1+K_{r\dz}\otimes Z_r$ for $r\geq 1$.
\end{itemize}
\end{Prop}

\begin{pf}By definition, it suffices to show that
each formula holds for any finite field $k$ with $q:=|k|\gg0$. Thus,
the formulas (2) and (3) follow directly from \cite{BurSch}.

We now prove (1). Since for $\ul{z}=\ul{z}_{\bfgz,q}$, both subsheaves and quotient sheaves of
$S_q(\bfgz,\ul{z})$ are again torsion sheaves with supports in
$\ul{z}$, we have
$$\aligned
{}& \Dz_q([S_q(\bfgz,\ul{z})]K_{\bfa})\\
=&\sum_{\bfaz,\bfbz}v_q^{\langle\bfaz,\bfbz\rangle} F_{S_q(\bfaz,\ul{z}),S_q(\bfbz,\ul{z})}^{S_q(\bfgz,\ul{z})}
\frac{a_{S_q(\bfaz,\ul{z})} a_{S_q(\bfbz,\ul{z})}}{a_{S_q(\bfgz,\ul{z})}}[S_q(\bfaz,\ul{z})]K_{\bfbz+\bfa}\otimes
[S_q(\bfbz,\ul{z})]K_{\bfa}\\
=&\sum_{\bfaz,\bfbz}v_q^{\langle\bfaz,\bfbz\rangle}\vphi_{\bfaz,\bfbz}^{\bfgz}(q)
\frac{a_{\bfaz}(q) a_{\bfbz}(q)}{a_{\bfgz}(q)}[S_q(\bfaz,\ul{z})]K_{\bfbz+\bfa}\otimes [S_q(\bfbz,\ul{z})]K_{\bfa}.
\endaligned$$
\end{pf}

\begin{Prop} The comultiplication $\Dz$ induces a topological bialgebra structure on the extended generic Hall algebra
$\wt\bfcalH_\bfv(\bbX)$.
\end{Prop}

\begin{pf} %{\blue Recall that the family $(u_{\bfaz}Z_{\ul{l}}K_{\bfbz})_{(\bfaz, \ul{l},\bfbz)\in\cS\times\mathcal{L}\times\cS}$
%is a $\bbQ[\bfv,\bfv^{-1}]$-basis of $\wt\calH_\bfv(\bbX)$. }
By Proposition \ref{comultiplication formula}, for $r\geq 1$ and $\bfaz\in\chi_t$, $\Dz(Z_r)$
and $\Dz(u_{\bfaz}K_{\bfa})$ belong to the completion
$\wt\bfcalH_\bfv(\bbX)\widehat{\otimes}\wt\bfcalH_v(\bbX)$. It suffices to show that
$\Dz(u_{\bfaz})\in\wt\bfcalH_\bfv(\bbX)\widehat{\otimes}\wt\bfcalH_\bfv(\bbX)$ for $\bfaz\in\chi_f$. Note
that for each $q\in\frakQ$, the category $\vect\bbX_q$ can be
generated by line bundles. It is reduced to prove that $\Dz(u_{[\co(\vx)]})$ belongs to
$\wt\calH_\bfv(\bbX)\widehat{\otimes}\wt\bfcalH_\bfv(\bbX)$. Since $\wt\bfcalH_\bfv(\bbX)$ is
closed under grading shift by $\vx$; see \cite[Cor.~5.18]{BurSch}, the assertion follows from
Propositions \ref{common subalgebra}, Proposition \ref{comultiplication formula}(2), and the fact that
$\Thz_{\vx}\in\bfcalH_\bfv(\bbX)$ for all $\vx\in\bbL_+$.
\end{pf}

\subsection{Drinfeld double of extended generic Hall algebra}
By \cite{Green}, for each $q\in\frakQ$, there
is a paring (called the Green's pairing)
$$\{-,-\}_q: \wtH(\bbX_q)\times \wtH(\bbX_q)\lra \bbQ[v_q, v_q^{-1}]$$
on $\wtH(\bbX_q)$ defined by
$$\{[S]K_{\bfa},[S']K_{\bfb}\}_k=v_q^{(\bfa,\bfb)}\frac{\delta_{S,S'}}{a_{S}}\;\;
\forall\,S,S'\in\coh\bbX_q,\;\bfa,\bfb\in K_0(\bbX).$$
Moreover, this pairing is non-degenerate and symmetric, and satisfies
$$\{a b,c\}_q=\{a\otimes b,\Dz_q(c)\}_q\;\,\forall\,a,b,c\in \wtH(\bbX_q).$$
 They give rise to a pairing
$$\{-,-\}: \prod_q\wtH(\bbX_q)\times\prod_q\wtH(\bbX_q)\lra \prod_q\bbQ[v_q, v_q^{-1}]$$
 defined by
$$\{(a_q)_q,(b_q)_q\}=(\{a_q,b_q\}_q)_q.$$
 It is easy to see that $\bigoplus_q\bbQ[v_q ,v_q^{-1}]$ is an ideal of
$\prod_q\bbQ[v_q, v_q^{-1}]$, and we denote the
quotient ring by $\tbbQ$. Finally, we obtain a pairing
\begin{equation} \label{pairing-infinite-product}
\{-,-\}:\big(\prod_q\wtH(\bbX_q)/\wt\cI\big)\otimes \big(\prod_q\wtH(\bbX_q)/\wt\cI\big)\lra \tbbQ
\end{equation}
satisfying that
$$\{(a_q)_q (b_q)_q,(c_q)_q\}=\{(a_q)_q\otimes (b_q)_q,\Dz((c_q)_q)\}_q,\;\forall\,(a_q)_q,(b_q)_q,(c_q)_q\in\prod_q\wtH(\bbX_q)/\wt\cI.$$

Now we recall some notations from \cite{BurSch}. For each $1\leq i\leq t$ and $r>0$,
define
$${\rm def}_{x_i}=\frac{[r]^2}{r}\big(\frac{1}{\bfv^{2rp_i}-1}-\frac{1}{\bfv^{2r}-1}\big)$$
which is the defect of the exceptional point $x_i$ when specializing $\bfv$ to $v_q$, where
$[r]=(\bfv^r-\bfv^{-r})/(\bfv-\bfv^{-1})$.

\begin{Prop} \label{pairing-value-1}
We have the following equalities:
\begin{itemize}
\item[(1)] $\{Z_r, Z_s\}=\dz_{r,s}(\frac{1}{\bfv-\bfv^{-1}}\frac{[2r]}{r}+\sum_{i\in I}{\rm def}_{x_i})$, $\forall\,r,s\geq 1$;
\item[(2)] $\{u_{\bfaz}K_{\bfa}, u_{\bfbz}K_{\bfb}\}=\bfv^{(\bfa,\bfb)}\frac{\delta_{\bfaz,\bfbz}}{a_{\bfaz}}$,
           $\forall\,\bfaz,\bfbz\in\cS,\;\bfa,\bfb\in K_0(\bbX)$;
\item[(3)] For $\vx=\sum_{i\in I}l_i\vx_i+l\vc\in\bbL_+$ and
$\bfaz=(\alpha, \lz)\in\cS$ with $\lz=((\lz^{(1)}, d_1), \ldots, (\lz^{(r)},d_r))$,
$\{u_{\bfaz}, \Thz_{\vx}\}=0$ or equals to
$$\frac{\bfv^{(\bfaz,\bfaz)+l+m}}{a_{\bfaz}}\prod_{1\leq s\leq r,\lz^{(s)}\neq \emptyset}(1-\bfv^{-2d_s})
\prod_{i\in I,(m_i,l_i)\neq (0,0)}(1-\bfv^{-2});$$
 see \eqref{def-Theta_x} for the notations.
\end{itemize}
\end{Prop}

\begin{pf} By definition, it suffices to show that
each equality holds for $q\in\frakQ$ with $q\gg0$. Then
(1) follows from \cite[Prop.~6.3]{BurSch} and (2) can be proved by an
easy calculation.

(3) For each $q\in\frakQ$ and $\ul{z}=\ul{z}_{\bfaz,q}\in\cX_{\bbF_q}(\ul{d})$,
$$S_q(\bfaz, \ul{z})=S_q(\az)\oplus(\oplus_{1\leq s\leq r}S_q(\lz^{(s)}, z_s)),$$
so $\{[S_q(\bfaz, \ul{z})],\Thz_{\vx,q}\}_q\neq 0$ if and only if $[S_q(\bfaz, \ul{z})]$ occurs
in the expression of $\Thz_{\vx,q}$ in \eqref{def-Theta_x}. This happens only if
$\det\bfaz:=\det\alpha+\sum_{1\leq s\leq r}|\lz^{(s)}|d_s\vc=\vx$ and
the length of $\lz^{(s)}$ is $1$ for each $1\leq s\leq r$.
In this case, the coefficient of $[S_q(\bfaz, \ul{z})]$ is given by
$$v_q^{l+m}\prod_{1\leq s\leq r,\lz^{(s)}\neq \emptyset}(1-v_q^{-2d_s})
\prod\limits_{i\in I,(m_i,l_i)\neq (0,0)}(1-v_q^{-2}).$$
On the other hand,
$$\{[S_q(\bfaz, \ul{z})], [S_q(\bfaz, \ul{z})]\}_q=\frac{v_q^{(\bfaz,\bfaz)}}{a_{S_q(\bfaz, \ul{z})}}.$$
This completes the proof.

\end{pf}

Our next goal is to introduce the reduced Drinfeld double of the
topological bialgebra $\wt\bfcalH_\bfv(\bbX)=\wt\calH_\bfv(\bbX)\otimes_{\bbQ(\bfv)}\bbQ(\bfv)$.
 By Proposition \ref{pairing-value-1}, the pairing in \eqref{pairing-infinite-product}
 induces a pairing
$$\{-,-\}: \wt\bfcalH_\bfv(\bbX)\times \wt\bfcalH_\bfv(\bbX)\lra \bbQ(\bfv).$$

Consider the pair of $\bbQ(\bfv)$-algebras $\wt\bfcalH^{\pm}_\bfv(\bbX)$ with bases
$$\big\{u^{\pm}_{\bfaz}Z^{\pm}_{\ul{l}}K^{\pm}_{\bfa}\mid \bfaz\in\cS, \ul{l}\in{\cal L}, \bfa\in K_0(\bbX)\big\}.$$
 The \emph{Drinfeld double} of
the topological bialgebra $\wt\bfcalH_\bfv^\pm(\bbX)$ with respect to the pairing
$\{-,-\}$ is the associative algebra $\wt{D}\bfcalH_\bfv(\bbX)$, defined
as the free product of algebras $\wt\bfcalH^+_\bfv(\bbX)$ and $\wt\bfcalH^-_\bfv(\bbX)$
subject to the relations
$$R(a,b):\;\;\sum a_1^- b_2^+\{a_2,b_1\}=\sum b_1^+a_2^-\{a_1,b_2\},$$
 where $a,b\in \wt\bfcalH_\bfv(\bbX)$, $\Delta(a)=\sum a_1\otimes a_2$, and $\Delta(b)=\sum b_1\otimes b_2$.
The \emph{reduced Drinfeld double} $D\bfcalH_\bfv(\bbX)$ is the quotient of
$\wt{D}\bfcalH_\bfv(\bbX)=\wt\bfcalH^+_\bfv(\bbX)\otimes\wt\bfcalH^-_\bfv(\bbX)$ by the Hopf ideal generated by
$K_{\bfa}^+\otimes K_{-\bfa}^--1\otimes 1$ for $\bfa\in K_0(\bbX)$. We then have the following triangular
decomposition of $D\bfcalH_\bfv(\bbX)$ as $\bbQ(\bfv)$-vector spaces
$$\bfcalH^+_\bfv(\bbX)\otimes_{\bbQ(\bfv)}\bfcalK\otimes_{\bbQ(\bfv)}\bfcalH^-_\bfv(\bbX)\stackrel{\cong}{\lra} D\bfcalH_\bfv(\bbX),$$
where $\bfcalK=\bbQ(\bfv)[K_0(\bbX)]$. For each $\bfa\in K_0(\bbX)$, we will write
$$K_\bfa=K_\bfa^+\otimes 1=1\otimes K_{\bfa}^-\in D\bfcalH_\bfv(\bbX).$$

\begin{Rems} (1) Let $\bfcalC_\bfv(\bbX)$ denote the $\bbQ(\bfv)$-subalgebra of
$\bfcalH_\bfv(\bbX)=\calH_\bfv(\bbX)\otimes\bbQ(\bfv)$
generated by $u_{[\co(l\vc)]}$ for $l\in\bbZ$,
$u_{[S_{i,j}]}$ for $i\in I$, $0\leq j\leq p_i-1$, and $T_r$ for
$r\geq 1$, called the {\it generic composition algebra} of $\bbX$. It was shown that
\cite[Th.~5.3]{Sch2004} that $\bfcalC_\bfv(\bbX)$ is isomorphic to
a quantized loop algebra ${\bf U}_{\bfv}(\widehat{\fn})$.

(2) In \cite{DJX}, Dou, Jiang and Xiao proved that for any field $k$ with
$q$ elements, the double Ringel--Hall algebra of $H(\bbX_k)$ provides a realization
of the quantized loop algebra ${\bf U}_{v_q}(\sL\fg)$ (over $\bbC$, specialized at
$\bfv=v_q=\sqrt{q}$) of the simple Lie algebra $\fg$ associated with $\bbX_k$.
Let $D\bfcalC_\bfv(\bbX)$ be the subalgebra of $D\bfcalH_\bfv(\bbX)$ generated by the set
$$\big\{u_{[\co(l\vc)]}^{\pm},u_{[S_{i,j}]}^{\pm}, T_{r}^{\pm}, K_{\bfa}^{\pm} \mid
l\in\bbZ, i\in I, 0\leq j\leq p_i-1, r\geq 1, \bfa\in K_0(\bbX)\big\},$$
called the \emph{generic Drinfeld double composition algebra} of $\bbX$.
We would expect a generic version of the isomorphism given in \cite[Th.~5.5]{DJX}, that is, there
is a $\bbQ(\bfv)$-algebra isomorphism ${\bf U}_{\bfv}(\sL\fg)\to D\bfcalC_\bfv(\bbX)$.

\end{Rems}

\section{Hall polynomials for tame quivers}

In this section we first define decomposition sequences for a tame quiver $Q$ analogously as in Section 3
and then use Theorem \ref{existence-HallPoly} together with \cite[Prop.~5]{cram}
to prove the existence of Hall polynomials for each triple of decomposition
sequences for $Q$. This not only refines the main theorem in \cite{Hub}, but also
confirms a conjecture in \cite[Conj.~3.4]{BG}.

Throughout this section, $Q=(Q_0,Q_1)$ denotes an acyclic tame quiver, that is, $Q$ contains
no oriented cycles and the underlying diagram of $Q$ is an extended Dynkin diagram of type $\widetilde A,
\widetilde D$ and $\widetilde E$. Let $kQ$ be the path algebra of $Q$ over a field $k$.
By ${\rm mod}\,kQ$ we denote the category of finite dimensional (left) $kQ$-modules.
It is well known from \cite{DR76} that the subcategory $\text{ind\,}kQ$ of indecomposable
$kQ$-modules admits a disjoint decomposition
$$\text{ind\,}kQ=\cP_k\cup\cR_k\cup\cI_k,$$
 where $\cP_k$ (resp., $\cR_k,\cI_k$) denotes the subcategories of indecomposable preprojective
(resp., regular, preinjective) modules. Moreover, $\cR_k$ consists of finitely many non-homogeneous
tubes and infinitely many homogeneous tubes which are parameterized by a subset $\bbH_k$
of $\bbP^1(k)$.

Similar to Section 3, we denote by $\Xi=\Xi(kQ)$ the
set of isoclasses of $kQ$-modules which clearly depends on the
ground field $k$. Let $\Xi_{\cP}$, $\Xi_{\cR}$ and $\Xi_{\cI}$ be
the subsets of $\Xi$ consisting of the isoclasses of preprojective,
regular and preinjective modules, respectively. Further, let
$\Xi_\nh$ be the subset of $\Xi$ formed by the isoclasses of $kQ$-modules without
homogeneous regular summands. Hence, the set $\Xi_\nh$ can be
described combinatorially and is independent of $k$. Moreover, each
module whose summands lie in a single homogeneous tube is
determined by a partition. For each $\az\in\Xi$, we fix a
representative $M_k(\az)$ in the class $\az$. Given
$\az,\bz\in\Xi$, we write $\az\oplus\bz$ for the isoclass of
$M_k(\az)\oplus M_k(\bz)$. Thus, each $\az\in\Xi$ can be uniquely
written as $\az=\az_{\cP}\oplus\az_{\cR}\oplus\az_{\cI}$ with
$\az_{\cP}\in\Xi_{\cP}$, $\az_{\cR}\in\Xi_{\cR}$, and $\az_{\cI}\in\Xi_{\cI}$.

A {\it decomposition sequence of type $\ul{d}=(d_1,\ldots,d_r)$} for $Q$ is a
pair $\bfaz=(\alpha, \lz)$, where $\alpha\in\Xi_\nh$ and $\lz$ is a
Segre sequence of type $\ul{d}$ (see Section 3 for the definition).
For each Segre sequence $\lz$ of type $\ul{d}$ and $\ul{z}=(z_1,
\ldots, z_r)\in {\cal X}_k(\ul{d})$, define
$$M_k(\lz, \ul{z}):=\bigoplus\limits_{i=1}^{r}M_k(\lz^{(i)}, z_i)\in{\rm mod}\,kQ,$$
where $M_k(\lz^{(i)}, z_i)$ is the regular $kQ$-module in the homogeneous tube
associated with the point $z_i$ determined by the partition
$\lz^{(i)}$. Further, for each decomposition sequence
$\bfaz=(\alpha, \lz)$ of type $\ul{d}$, define
$$M_k(\bfaz, \ul{z}):=M_k(\alpha)\oplus M_k(\lz, \ul{z}).$$
 Finally, by $\cM=\cM(Q)$ we denote the set of all decomposition sequences for $Q$
via identifying $(\az,\lz)$ and $(\az,\mu)$, where $\lz$ is obtained from $\mu$
by inserting and removing some pairs $(\emptyset, d)$; see Remark \ref{equiv-Segre-class}.
Further, we define the subsets $\cM_{\cP}$, $\cM_{\cR}$ and $\cM_{\cI}$ of $\cM$ in a natural sense.
Each $\bfaz\in\cM$ can be written as $\bfaz=\bfaz_{\cP}\oplus\bfaz_{\cR}\oplus\bfaz_{\cI}$
with $\bfaz_{\cP}\in\cM_{\cP},$ $\bfaz_{\cR}\in\cM_{\cR}$ and $\bfaz_{\cI}\in\cM_{\cI}$.

\begin{Defn} \label{Def-Hall-poly for tame quiver} Given $\bfaz, \bfbz,\bfgz\in\cM=\cM(Q)$ of type $\ul{d}$,
if there exists a polynomial $\psi_{\bfaz,\bfbz}^{\bfgz}\in\bbZ[T]$
such that for each finite field $k$ with $q:=|k|\gg0$,
$$\psi_{\bfaz,\bfbz}^{\bfgz}(q)= F_{M_k(\bfaz, \ul{z}), M_k(\bfbz, \ul{z})}^{M_k(\bfgz, \ul{z})}
\quad\textrm{for all }\ul{z}\in \cX_k(\ul{d}),$$
 then the Hall polynomial $\psi_{\bfaz,\bfbz}^{\bfgz}$ is  said to exist for
$\bfaz$, $\bfbz$ and $\bfgz$.
\end{Defn}

The main aim of this section is to prove the existence of Hall polynomials for
the tame quiver $Q$. We need some preparation.

Let $\bbX_k$ be a weighted projective line associated with
the tame quiver $Q$ in the sense that there exists an equivalence
$D^b(\coh\bbX_k)\cong D^b(kQ\text{-mod})$ of the bounded derived categories.
More precisely, there exists a tilting bundle $T_{k}$ in $\coh\bbX_k$, whose
summands form a complete slice in the Auslander--Reiten quiver of
$\vect\bbX_k$, such that the endomorphism algebra of $T_k$ is
isomorphic to the path algebra $kQ$.
Let $\sF(T_k)$ be the torsion-free class of $\coh\bbX_k$ induced by
$T_k$ which consists of the sheaves $S_k\in\coh\bbX_k$
satisfying $\Hom(T_k, S_k)=0$. Since $T_k$ is a vector
bundle, it is easily seen that $\sF(T_k)$ is a full subcategory of
$\vect\bbX_k$, which can be described combinatorially and is independent of the field $k$. Thus, we will
drop the index $k$ and view $\sF(T)$ as a subset of $\cS$. Let ${\scr T}(T)$ be the
subset of $\cS$ consisting of $\bfaz\in \cS$ such that $M_k(\bfaz, \ul{z})$ is generated by $T_k$
for each $\ul{z}\in \cX_k(\ul{d})$, where $\ul{d}$ is the type of $\bfaz$. Then
each decomposition sequence $\bfaz$ admits a decomposition $\bfaz=\bfaz_+\oplus\bfaz_-$  with
$\bfaz_+\in {\scr T}(T)$ and $\bfaz_-\in\sF(T)$. Moreover, $T$ induces
bijections $\Psi_1: {\scr T}(T)\to\cM_{\cP\oplus\cR}$
(preserving the Segre sequence) and $\Psi_2:\sF(T)\to\cM_{\cI}$.
As a consequence, the functor $\Hom(T,-)$ induces a triangulated equivalence of bounded derived categories
$${\bf R}\Hom(T,-): D^b(\coh\bbX_k)\lra D^b(kQ\text{-mod}).$$

Let $D\bfH(kQ)$ and $D\bfH(\bbX_k)$ denote the double Ringel--Hall algebras of $kQ$ and $\bbX_k$
over the field $\bbQ(v_q)$, respectively; see \cite{X97,BurSch,DJX}.
Applying \cite[Prop.~5]{cram} gives the following result.

\begin{Lem} \label{Cramer-result} There exists a $\bbQ(v_q)$-algebra
isomorphism
$$D\bfH(kQ)\stackrel{\cong}{\lra}D\bfH(\bbX_k),$$
which takes
$$\aligned{}
&[M_k({\bfthz_{\cP}\oplus\bfthz_{\cR}},\ul{z})]^+ \lmto  [S_k({\bfthz_{+}},\ul{z})]^+,\\
&[M_k(\bfthz_{\cI},\ul{z})]^+ \lmto v_q^{-\langle \bfthz_{-},
\bfthz_{-}\rangle}[S_k(\bfthz_{-},\ul{z})]^-K_{\bfthz_{-}}\,\;;
\endaligned$$
where $\bfthz=\bfthz_{\cP}\oplus\bfthz_{\cR}\oplus\bfthz_{\cI}\in\cM$ of type $\ul{d}$, $\ul{z}\in \cX_k(\ul{d})$,
$\bfthz_+=\Psi_1^{-1}(\bfthz_{\cP}\oplus\bfthz_{\cR})$, and $\bfthz_-=\Psi_2^{-1}(\bfthz_{\cI})$.
\end{Lem}

Now let $\bfcalH_\bfv(T)$ be the $\bbQ(\bfv)$-submodule of $D\bfcalH_\bfv(\bbX)$ spanned by the set
$$\{u^+_{\bfaz}Z^+_{\ul{l}}u^-_{\bfbz}K_{\bfa}\mid  \bfaz\in\sT(T),\,\ul{l}\in{\cal L}, \, \bfbz\in \sF(T),\,
\bfa\in K_0(\bbX)\}.$$

\begin{Prop}
$\bfcalH_\bfv(T)$ is a $\bbQ(\bfv)$-subalgebra of $D\bfcalH_\bfv(\bbX)$.
\end{Prop}

\begin{pf} Since the torsion (resp., torsion-free) class induced by $T$ is
closed under extension, the $\bbQ(\bfv)$-submodule of $D\bfcalH_\bfv(\bbX)$
spanned by $u^+_{\bfaz}Z_{\ul{l}}^+$ for $\bfaz\in\sT(T)$
(resp., by $u^-_{\bfbz}K_{\bfa}$ for $\bfbz\in\sF(T)$ and $\bfa\in K_0(\bbX)$) is a subalgebra of $D\bfcalH_\bfv(\bbX)$.
Hence, we only need to show that $u_{\bfbz}^- u_{\bfaz}^+Z^+_{\ul{l}}\in\bfcalH_\bfv(T)$ for
$\bfaz\in\sT(T)$, $\ul{l}\in{\cal L}$, and $\bfbz\in\sF(T)$.

Take $\bfbz\in\sF(T)$ and write $\Dz(u_{\bfbz})=\sum b_1\otimes b_2$. Since $\sF(T)$ is closed
under subobjects, all the $b_2$ are generated by $u_{\bfbz}$ with $\bfbz\in\sF(T)$.
Then
$$u_{\bfbz}^- Z_r^{+}=Z_r^+ u_{\bfbz}^- + \sum\{Z_r, b_1\} b_2^- \in \bfcalH_\bfv(T).$$
 It remains to show $u_{\bfbz}^- u_{\bfaz}^+\in\bfcalH_\bfv(T)$ for
$\bfaz\in\sT(T),\,\bfbz\in\sF(T)$.
Assume $\Dz(u_{\bfaz})=\sum a_1\otimes a_2$. Then the $a_1$ are generated by
$u_{\bfaz},\bfaz\in\sT(T)$ and $Z_r, r\geq 1$.
By the defining relation $R(u_{\bfaz}, u_{\bfbz})$ in $D\bfcalH_v(\bbX)$, we have
\begin{equation}\label{drinfeld double}
u_{\bfbz}^- u_{\bfaz}^+=\sum \{a_2, b_1\} a_1^+ b_2^-,
\end{equation}
 where $\Dz(u_{\bfbz})=\sum b_1\otimes b_2$ and $a_1^+ b_2^-\in\bfcalH_v(T)$.
We need to show that the sum is finite. Since $u_{\bfaz}=u_{\bfaz_f}u_{\bfaz_t}$,
we only need to consider the following two cases:

\medskip

\emph{Case 1}: $\bfaz\in\cS_t$. In this case, $a_1$ and $a_2$ are generated by
$u_{\bfaz}$, $\bfaz\in\sT(T)$. Hence, there are only finitely
many choices of $a_1$ and $a_2$ since their degrees are bounded by $\deg\bfaz$.
Furthermore, there are only finitely many $b_1$'s such that $\{a_2,
b_1\}\neq 0$. This forces that there are finitely many choices of $b_2$ which give nonzero
terms in the right hand side of \eqref{drinfeld double}.

\medskip

\emph{Case 2}: $\bfaz\in\cS_f$. In this case, each term $a_2$ can be assumed to have the
form $u_{\bfgz}$ for some $\bfgz\in\cS_f$. Then $\{a_2, b_1\}\neq 0$ implies that $b_1=c_{\bfgz}u_{\bfgz}$
for some $c_{\bfgz}\in\bbQ(\bfv)$. Thus, there are an epimorphism
$S_k(\bfbz)\twoheadrightarrow S_k(\bfgz)$ and a monomorphism $S_k(\bfgz)\rightarrowtail S_k(\bfaz)$,
which ensures that there are only finitely
many choices of $\bfgz$ and so for $a_2$. Hence, there are
finitely many triples $(a_2,b_1,b_2)$ which contribute a nonzero term in \eqref{drinfeld double}.
\end{pf}

Since all the exceptional simple sheaves belong to the torsion classes, we obtain by an argument similar
to that in the proof of Proposition \ref{three-bases-generic-HA} that the two sets
$$\{u^+_{\bfaz}T^+_{\ul{l}}u^-_{\bfbz}K_{\bfa}\}\;\text{ and }\;\{u^+_{\bfaz}\Thz^+_{\ul{l}}u^-_{\bfbz}K_{\bfa}\},$$
where $\bfaz\in\sT(T),\,\ul{l}\in {\cal L}, \, \bfbz\in\sF(T)$, and
$\bfa\in K_0(\bbX)$, are both $\bbQ(\bfv)$-bases of $\bfcalH_\bfv(T)$.

A rational function $\phi(\bfv)=f(\bfv)/g(\bfv)\in\bbQ(\bfv)$ is said to be {\it nearly integral}
if $f(\bfv),g(\bfv)\in\bbZ[\bfv]$ and $g(\bfv)$ is monic. We say that an element $a\in\bfcalH_\bfv(\bbX)$
is generated by a set $X$ with nearly integral coefficients if $a$ is a linear combination
of monomials of elements in $X$ whose coefficients are nearly integral functions.

\begin{Lem}\label{vb vs line bds with int/monic} For each $\bfthz\in\sF(T)$,
the element $u_{\bfthz}\in \bfcalH_\bfv(\bbX)$ can be generated by
$\{u_{[\co(\vx)]}\mid \vx\in\bbL\}$ with nearly integral coefficients.
\end{Lem}

\begin{pf} By the assumption, $S_k(\bfthz)=E$ is a vector bundle. If $E$ is
decomposable, then $E$ can be written as $E=\bigoplus\limits_{1\leq i\leq r}E_i^{l_i}$, where $E_i$ are indecomposable
satisfying $\Hom(E_i,E_j)=0$ for $i>j$. Then
$$u_{[E]}=\bfv^{-\sum\limits_{1< i\leq r}l_i l_{i-1}-\sum\limits_{i<j} l_i l_j\langle E_i, E_j\rangle}
\frac{u_{[E_1]}^{l_1}}{[l_1]!}\ldots \frac{u_{[E_r]}^{l_r}}{[l_r]!}.$$
Thus, it suffices to prove the assertion in the case where $E$ is indecomposable.

Now suppose that $E$ is indecomposable. By Lemma \ref{vector bundle factorization},
for each finite field $k=\bbF_q$, there is an exact sequence
$$0\lra L\lra E\lra F\lra 0$$
in $\vect\bbX_k$ such that $L$ is a line bundle and $\Ext^1(F, L)\cong k$. Then,
in the Ringel--Hall algebra $H(\bbX_k)$, we have the equalities
$$[L][F]=v_q^{\langle L, F\rangle}[L\oplus F]\;\text{ and}$$
$$ [F][L]=v_q^{\langle F, L\rangle}(v_q^{2\langle L, F\rangle} [L\oplus F]+\frac{v_q^2-1}{a_F}[E])
=v_q^{-1}(v^{\langle L, F\rangle} [L][F]+\frac{v_q^2-1}{a_F}[E]).$$
It follows that
$$[E]=\frac{v_q}{v_q^2-1} a_F[F][L]-\frac{v_q^{\langle L, F\rangle}}{v_q^2-1} a_F[L][F].$$
In other words, we have in $\bfcalH_\bfv(\bbX)$,
$$u_{[E]}=\frac{\bfv}{\bfv^2-1} a_{\bfthz'}u_{[F]}u_{[L]}-\frac{\bfv^{\langle L, F\rangle}}{\bfv^2-1} a_{\bfthz'} u_{[L]}u_{[F]},$$
 where $\bfthz'\in \sF(T)$ is defined by $S(\bfthz')=F$.
Since $F$ is a vector bundle with rank smaller than that of $E$, the assertion follows from an induction on the rank
of $E$.
\end{pf}

\begin{Lem}\label{commutative with int/monic}
For decomposition sequences $\bfthz\in\sF(T)$ and $\bfdz\in\sT(T)$, the element
$u^-_{\bfthz} u^+_{\bfdz}$ has nearly integral coefficients with respect to the basis
\begin{equation}\label{basis-H(T)}
\{u^+_{\bfaz}\Thz^+_{\ul{l}}u^-_{\bfbz}K_{\bfa}\mid  \bfaz\in\sT(T),\, \ul{l}\in{\cal L}, \,\bfbz\in\sF(T),\, \bfa\in K_0(\bbX)\}.
\end{equation}
\end{Lem}

\begin{pf} By Lemma \ref{vb vs line bds with int/monic}, for each $\bfgz\in\sF(T)$, $u_{\bfgz}$
can be generated by $\{u_{[\co(\vx)]}\mid \vx\in \bbL\}$ with nearly integral coefficients.
Then by Proposition \ref{comultiplication formula},
\begin{equation}\label{coproduct of vector bundles}
\Delta(u_{\bfgz})=\sum c_{\bfgz} a_{1,\bfgz}\otimes a_{2,\bfgz}
 \end{equation}
where each $a_{1,\bfgz}$ can be generated by $\{u_{[\co(\vx)]}, \Thz_{\vec y} \mid \vx\in \bbL, {\vec y}\in\bbL_+\}$
with nearly integral coefficients, and each $a_{2,\bfgz}$ has the form $u_{\bfgz'}$
for some $\bfgz'\in\sF(T)$. Since $u_{[\co(\vx)]}$ (resp.,  $\Thz_{\vx})$ can be
generated by $u_{[\co(l\vc)]}$'s (resp., $\Thz_{r\vc}$'s) and $u_{[S_{i,j}]}$'s
with nearly integral coefficients, we conclude that the $a_{1,\bfgz}$ in
 \eqref{coproduct of vector bundles} can be generated by
$$\{u_{[\co(l\vc)]}, u_{[S_{i,j}]}, \Thz_{r\vc} \mid l\in \bbZ, r\in\bbN, i\in I, 0\leq j\leq p_i-1\}$$
with nearly integral coefficients.

We first consider the following two special cases of $\bfdz$.

\emph{Case 1}: $\bfdz\in\cS_t$. In this case, $\Delta(u_{\bfdz})=\sum c_{\bfdz} u_{\bfdz_1}\otimes u_{\bfdz_2}$,
 where $\bfdz_1,\bfdz_2\in\cS_t$ and $c_{\bfdz}=\bfv^{\langle u_{\bfdz_1},u_{\bfdz_2}\rangle}F_{\bfdz_1,\bfdz_2}^{\bfdz}$
 are nearly integral functions. Hence,
 $$u^-_{\bfthz} u^+_{\bfdz}=\sum c_{\bfthz}c_{\bfdz}\{a_{1,\bfthz}, u_{\bfdz_2}\} u^+_{\bfdz_1} a_{2,\bfthz}^-,$$
where the coefficients $c_{\bfthz}c_{\bfdz}\{a_{1,\bfthz}, u_{\bfdz_2}\}$ are nearly integral
functions by Proposition \ref{pairing-value-1}(3).

\emph{Case 2}: $\bfdz\in\cS_f$. In this case,
$u^-_{\bfthz} u^+_{\bfdz}=\sum c_{\bfthz}c_{\bfdz}\{a_{1,\bfthz}, a_{2,\bfdz}\} a^+_{1, \bfdz} a_{2,\bfthz}^-$.
Note that $a_{2,\bfdz}$ has the form $u_{\bfdz'}$ for some $\bfdz'\in\cS_f$. Thus,
$\{a_{1,\bfthz}, a_{2,\bfdz}\}\neq 0$ if and only if $u_{\bfdz'}$ is a nonzero term in the linear
combination of $a_{1,\bfthz}$ with respect to the third basis in Proposition \ref{three-bases-generic-HA}.
By Proposition \ref{pairing-value-1}, the coefficients
$c_{\bfgz}c_{\bfdz}\{a_{1,\bfthz}, a_{2,\bfdz}\}$ are nearly integral functions.

\medskip

In general, $u_{\bfdz}$ has a decomposition $u_{\bfdz}=u_{\bfdz_f}u_{\bfdz_t}$, where $\bfdz_f\in\cS_f$
and $\bfdz_t\in\cS_t$. Hence,
$$u^-_{\bfthz} u^+_{\bfdz}=u^-_{\bfthz}u^+_{\bfdz_f}u^+_{\bfdz_t}
=\sum c_{\bfthz}c_{\bfdz_f}\{a_{1,\bfthz}, a_{2,\bfdz_f}\} a^+_{1, \bfdz_f} a_{2,\bfthz}^-u^+_{\bfdz_t}$$
$$=\sum c_{\bfthz}c_{\bfdz_f} c_{\bfthz'}c_{\bfdz_t}\{a_{1,\bfthz}, a_{2,\bfdz_f}\}\{a_{1,\bfthz'}, u_{\bfdz_{t,2}}\} a^+_{1, \bfdz_f}u^+_{\bfdz_{t,1}} a_{2,\bfthz'}^-,$$
where each $a_{2,\bfthz}$ takes the form $u_{\bfthz'}$ for some $\bfthz'\in\sF(T)$,
and the coefficients are nearly integral.

Finally, by the construction of $T_r$ given in \cite[Sect.~6]{BurSch}, we obtain that for each $i\in I$,
$[u_{[S_{i,j}]}, T_r]$ can be generated by $u_{[S_{i,s}]}, 0\leq s\leq p_i-1$, with nearly
integral coefficients. An analogous result holds for $[u_{S_{i,j}}, \Thz_{r\vc}]$ by
using \eqref{relation of Tr and Thz_r}. This together with Proposition \ref{formula of product} (5)
completes the proof.

\end{pf}

\begin{Thm} \label{existence-HallPoly for tame quiver} For arbitrary $\bfaz, \bfbz,\bfgz\in\cM$
of type $\ul{d}$, the Hall polynomial $\psi_{\bfaz,\bfbz}^{\bfgz}$ exists.
\end{Thm}

\begin{pf} For each $\bfthz\in\cM$ of type $\ul{d}$, write
$\bfthz=\bfthz_{\cP}\oplus\bfthz_{\cR}\oplus\bfthz_{\cI}$. Further, set
$$\bfthz_+=\Psi_1^{-1}(\bfthz_{\cP}\oplus\bfthz_{\cR})\in\sT(T)\;\text{ and }\;
\bfthz_-=\Psi_2^{-1}(\bfthz_{\cI})\in\sF(T).$$
%Then by definition, $u^+_{\bfthz_{+}}u^-_{\bfthz_{-}}\in\bfcalH_v(T)$.
%For arbitrary $\bfaz, \bfbz,\bfgz\in\cD$ of type $\ul{d}$,
Combining with Proposition \ref{formula of product} and Lemma \ref{commutative with int/monic},
we obtain that in the expression of the product
$(u^+_{\bfaz_{+}}u^-_{\bfaz_{-}})(u^+_{\bfbz_{+}}u^-_{\bfbz_{-}})$
with respect of the basis \eqref{basis-H(T)} of $\bfcalH_{\bfv}(T)$, the coefficient of
$u^+_{\bfgz_{+}}u^-_{\bfgz_{-}}K_{\bfaz_{-}+\bfbz_{-}-\bfgz_{-}}$
is a nearly integral function $\xi_{\bfaz, \bfbz}^{\bfgz}(\bfv)$. For each finite field $k$ with $q=|k|\gg0$
and $\ul{z}\in \cX_k(\ul{d})$, we take a total ordering $\preccurlyeq$ in $\bbH_k$ so that
$\ul{z}=\ul{z}_{\bfaz,q}=\ul{z}_{\bfbz,q}=\ul{z}_{\bfgz,q}$ in defining the embedding $\Phi$ in
Proposition \ref{embedding-Phi}. Then $\xi_{\bfaz, \bfbz}^{\bfgz}(v_q)$ is the coefficient of
the basis element $[S_k({\bfgz_{+}},\ul{z})]^+[S_k({\bfgz_{-}},\ul{z})]^-K_{\bfaz_{-}+\bfbz_{-}-\bfgz_{-}}$
in the expression of the product
$$([S_k({\bfaz_{+}},\ul{z})]^+[S_k({\bfaz_{-}},\ul{z})]^-)([S_k({\bfbz_{+}},\ul{z})]^+[S_k({\bfbz_{-}},\ul{z})]^-).$$
By Lemma \ref{Cramer-result}, there exists an integer $l(\bfaz,\bfbz,\bfgz)$, depending on $\bfaz,\bfbz$
and $\bfgz$, such that $v_q^{l(\bfaz,\bfbz,\bfgz)}\xi_{\bfaz, \bfbz}^{\bfgz}(v_q)$
is the coefficient of the basis element $[M_k({\bfgz},\ul{z})]^+$ in the product
$[M_k({\bfaz},\ul{z})]^+[M_k({\bfbz},\ul{z})]^+$. Set $\eta_{\bfaz,\bfbz}^{\bfgz}(\bfv)=\bfv^{l(\bfaz,\bfbz,\bfgz)}\xi_{\bfaz,
\bfbz}^{\bfgz}(\bfv)$, which is again a nearly integral function. Then, by the definition of the multiplication
in the Ringel--Hall algebra $H(kQ)$ of $kQ$, $\eta_{\bfaz, \bfbz}^{\bfgz}(\bfv)$ takes integer
values at $v_q=\sqrt{q}$ for prime powers $q\gg 0$. This forces that
$\eta_{\bfaz, \bfbz}^{\bfgz}(\bfv)$ is an integer polynomial in $\bfv^2$, that is, there is a polynomial
$\psi_{\bfaz, \bfbz}^{\bfgz}(T)\in\bbZ[T]$ such that
$\psi_{\bfaz, \bfbz}^{\bfgz}(\bfv^2)=\eta_{\bfaz, \bfbz}^{\bfgz}(\bfv)$, as desired.
\end{pf}

We now consider a special case of the theorem above. Let $\bfaz,\bfbz\in\cM$ be such that
$M_k(\bfaz)=I$ and $M_k(\bfbz)=P$ are preinjective and preprojective $kQ$-modules, respectively,
where $k$ is a finite field.
Further, let $d\geq 1$ and $\lz$ be a Segre sequence of type $(d)$ (i.e., $\lz$ is a partition).
Set $\bfgz=(0,\lz)\in\cM$. Then $\chi_k(d)$
consists of points in $\bbH_k$ of degree $d$ and for each $z\in \chi_k(d)$, $\cI_\lz(z):=M_k(\bfgz,z)$ is a
regular $kQ$-module whose summands all lie in the homogeneous tube corresponding to $z$. Applying
Theorem \ref{existence-HallPoly for tame quiver} gives the Hall polynomial $\psi^{\bfgz}_{\bfaz,\bfbz}$.
Therefore, for each finite field $k$ and $z_1,z_2\in\chi_k(d)$,
\begin{equation}\label{equality-HN-BG}
F^{\cI_\lz(z_1)}_{I,P}=F^{\cI_\lz(z_2)}_{I,P}.
\end{equation}
 Following the notation in \cite[2.5]{BG}, for three $kQ$-modules $A,B,C$, set
$$F^{A,B}_C=\frac{|\Ext^1_{kQ}(B,A)_C|}{|\Hom_{kQ}(B,A)|}.$$
 Since $|\Aut(\cI_\lz(z_1))|=|\Aut(\cI_\lz(z_2))|$ and $\Hom_{kQ}(I,P)=0$, applying
Lemma \ref{Riedtmann-Peng-formula} to \eqref{equality-HN-BG} gives the equality
$$F^{P,I}_{\cI_\lz(z_1)}=|\Ext^1_{kQ}(I,P)_{\cI_\lz(z_1)}|=|\Ext^1_{kQ}(I,P)_{\cI_\lz(z_2)}|=F^{P,I}_{\cI_\lz(z_2)}$$
 for all $z_1,z_2\in\chi_k(d)$. This is exactly the equality conjectured in \cite[Conj.~3.4]{BG}.

\begin{Cor} Conjecture {\rm 3.4} in \cite{BG} holds.
\end{Cor}

The following result is a consequence of Theorem \ref{existence-HallPoly for tame quiver} whose proof
is analogous to that of Corollary \ref{Hall-poly-3-obj}.

\begin{Cor} \label{Hall-poly-3-obj-tame-quiver} Let $k$ be a finite field and fix three $kQ$-modules $M,N,Z$.
Then there exists a polynomial $\psi_{M,N}^Z\in\bbZ[T]$ such that for each conservative field extension $K$
of $k$ relative to $\{M,N,Z\}$,
$$\psi_{M,N}^Z(|K|)=F_{M^K,N^K}^{Z^K}.$$
\end{Cor}

\vspace{.8cm}
 \noindent {\bf Acknowledgement.} We would like to thank Jie Xiao for stimulating discussions
 and helpful comments. Especially, the idea of studying Hall polynomials for
 tame quivers via weighted projective lines comes from his suggestion.

\end{document}